\documentclass[12pt,a4paper]{article}

\usepackage{amssymb}
\usepackage{amsmath}
\usepackage{amsthm}
\usepackage{xypic}

\theoremstyle{plain}
\newtheorem{theorem}{Theorem}
\newtheorem{proposition}[theorem]{Proposition}
\newtheorem{corollary}[theorem]{Corollary}
\theoremstyle{remark}
\newtheorem{remark}[theorem]{Remark}
\theoremstyle{definition}
\newtheorem{definition}[theorem]{Definition}

\newtheorem{lemma}[theorem]{Lemma}

\renewcommand{\baselinestretch}{1.2}

\def\varinjlim_#1{\lim\limits_{\longrightarrow\atop{#1}}}

\def\End{\mathop{\rm End}\nolimits}

\def\Aut{\mathop{\rm Aut}\nolimits}
\def\Hom{\mathop{\rm Hom}\nolimits}
\def\id{\mathop{\rm id}\nolimits}
\def\BSU{\mathop{\rm BSU}\nolimits}
\def\BU{\mathop{\rm BU}\nolimits}
\def\BPU{\mathop{\rm BPU}\nolimits}

\def\GBr{\mathop{\rm GBr}\nolimits}
\def\FB{\mathop{\rm FB}\nolimits}

\def\BFr{\mathop{\rm BFr}\nolimits}
\def\BN{\mathop{\rm BN}\nolimits}
\def\SU{\mathop{\rm SU}\nolimits}
\def\BSN{\mathop{\rm BSN}\nolimits}
\def\EU{\mathop{\rm EU}\nolimits}
\def\EN{\mathop{\rm EN}\nolimits}
\def\U{\mathop{\rm U}\nolimits}

\def\N{\mathop{\rm N}\nolimits}
\def\F{\mathop{\rm F}\nolimits}
\def\PGL{\mathop{\rm PGL}\nolimits}

\def\BGL{\mathop{\rm BGL}\nolimits}

\def\Fred{\mathop{\rm Fred}\nolimits}

\def\GL{\mathop{\rm GL}\nolimits}
\def\SL{\mathop{\rm SL}\nolimits}
\def\Gr{\mathop{\rm Gr}\nolimits}
\def\G{\mathop{\rm G}\nolimits}
\def\PU{\mathop{\rm PU}\nolimits}
\def\Z{\mathop{\rm Z}\nolimits}
\def\B{\mathop{\rm B}\nolimits}
\def\KSU{\mathop{\rm KSU}\nolimits}
\def\AB{\mathop{\rm AB}\nolimits}
\def\Br{\mathop{\rm Br}\nolimits}
\def\K{\mathop{\rm K}\nolimits}

\def\tr{\mathop{\rm tr}\nolimits}
\def\ker{\mathop{\rm ker}\nolimits}

\def\im{\mathop{\rm im}\nolimits}
\def\coker{\mathop{\rm coker}\nolimits}
\def\Fr{\mathop{\rm Fr}\nolimits}

\begin{document}

\author{A.V. Ershov}
\title{A generalization of the topological Brauer group}


\date{}
\maketitle
{\renewcommand{\baselinestretch}{1.0}

\begin{abstract}
In the present paper we study some homotopy invariants which can
be defined by means of bundles with fiber a matrix algebra. In particular,
we introduce some generalization of the Brauer group in the topological context
and show that any its element can be represented as a locally trivial bundle
with the structure group $\N_k^\times,\: k\in \mathbb{N}$ . Finally, we discuss its possible applications
in the twisted $K$-theory.
\end{abstract}}

\tableofcontents

\subsection*{Introduction}

The aim of this paper is to study different homotopy invariants
of $CW$-complexes which can be constructed by means
of locally-trivial bundles with fiber a matrix algebra.
An example of such invariant is the first obstruction for an
embedding of a locally-trivial matrix algebra bundle into a trivial one
(under some extra conditions which will be formulated below).
Another invariant we shall deal with is a homotopy functor $\GBr$
which can be treated as a generalization of the classical
``topological'' Brauer group. So let us remember
the basic idea of the Brauer group.

Suppose $X$ is a finite $CW$-complex.
The topological Brauer group $\Br(X)$
can be defined as the group of
equivalence classes of locally-trivial
bundles $A_k$ over $X$ with fibers $M_k(\mathbb{C})$
(for arbitrary $k\in \mathbb{N}$) with
respect to the following equivalence relation:
$$
A_k\sim B_l\; \Leftrightarrow \; A_k\otimes \End(\xi_m)\cong
B_l\otimes \End(\eta_n)
$$
for some complex vector bundles $\xi_m,\, \eta_n$
of rank $m,\, n$ respectively (in particular, $km=ln$).
Roughly speaking, we can say that $\Br(X)$ is
the group of obstructions for
the lifting of locally trivial bundles $A_k$ over $X$ with
fiber the matrix algebra $M_k(\mathbb{C})$ to bundles
of the form $\End(\xi_k)$ for some locally trivial
$\mathbb{C}^k$-bundle $\xi_k$ over $X.$

In order to define the generalized Brauer group $\GBr$
for any natural $k$ we introduce some $C^*$-algebra $\N_k\subset M_k({\cal B}(H))$,
where ${\cal B}(H)$ is the algebra of bounded operators in a separable complex Hilbert space $H$,
and consider bundles with the group of its invertible elements $\N_k^\times$ as a structure group.
Some of such bundles
come from finite-dimensional vector bundles due to a group homomorphism $\U(k)\rightarrow \N_k^\times$.
The generalized Brauer group is just the group of equivalence classes
of $\N_k^\times$-bundles modulo those that come from finite-dimensional vector bundles.

According to a classical theorem of J.-P. Serre,
there is an isomorphism
$\Br(X)\cong H^3_{tors}(X;\, \mathbb{Z})$
\cite{Grothendieck}.
In other words, all obstructions for
the lifting of locally-trivial bundles
$A_k$ over $X$ with
fiber $M_k(\mathbb{C})$ to bundles
of the form $\End(\xi_k)$ are independent of the
higher-dimensional integer cohomology of
dimensions greater than $3$.
In contrast to the classical case, the generalized
Brauer group $\GBr(X)$ actually depends on the
higher-dimensional cohomology $H^{2i+1}_{tors}(X;\, \mathbb{Z}),\:
i\geq 2$. On the other hand, just as in the classical group $\Br(X),$
any element of $\GBr(X)$ has finite order.

The last part of the paper is motivated by the recent progress in the twisted $K$-theory.
For any bundle with the structure group $\N_k^\times$ we define the corresponding twisted $K$-group.
We hope that our approach will provide an interesting
example of more general twistings than the ones considered up to now.

\smallskip

{\raggedright {\bf Acknowledgments}}\;
I am grateful to E.V. Troitsky
for constant attention to this work and all-round support
and to V.M. Manuilov and A.S. Mishchenko
for very helpful discussions.

\section{Matrix Grassmannians}
\label{appfl}

In this section we give definitions of objects we shall deal with below.
More precisely, we start with so-called ``matrix Grassmannians'' (which play
the same role for matrix algebras as the usual Grassmannians for vector spaces) and
try to develop a theory of matrix algebra bundles by analogy with the classical
theory of vector bundles. We restrict ourselves to considering a special class of matrix
algebra bundles (so-called ``floating'' bundles), because only in this case we avoid the localization
(hence a trivialization) of their theory. Finally in this way we obtain an invariant of $CW$-complexes
closely connected with the usual $K$-theory (cf. Theorem \ref{Kthcoincid}).

Note that in the course of the paper the basic field is $\mathbb{C}$.

\subsection{Basic definition}

Let us remember that the Grassmannian
$\G_{k,\, n},\, 0\leq k\leq n$ is a homogeneous space parameterizing
$k$-dimensional subspaces in the fixed
$n$-dimensional vector space $\mathbb{C}^n$.

By analogy with this definition we introduce a ``matrix Grassmannian''
which parameterizes central matrix subalgebras of a given dimension
in the fixed matrix algebra.

First, remember that a
{\it central} subalgebra is a subalgebra whose center coincides with
the center of the whole algebra (i.e. with the field of scalar matrices
in the case of matrix algebras). In particular, the identity matrix
of the ``big'' algebra is the identity element of any central subalgebra.
By $M_n(\mathbb{C})$ denote the algebra of all
$n\times n$ matrices over $\mathbb{C}.$

A central subalgebra in $M_n(\mathbb{C})$ isomorphic to $M_k(\mathbb{C})$
is called a $k$-{\it subalgebra}, for short. Note that such a subalgebra
exists only if $k|n,$ i.e. $n=kl$ for some natural $l$.

\begin{definition}
The {\it matrix Grassmannian}
$\Gr'_{k,\, l}$ is a homogeneous space
which parametrizes $k$-subalgebras in a fixed
matrix algebra $M_{kl}(\mathbb{C})$.
\end{definition}

\begin{proposition}
For any pair $k,\, l>1$ there exists the
matrix Grassmannian $\Gr_{k,\, l}'$.
\end{proposition}
{\noindent \it Proof}\; follows from Noether-Skolem's theorem
\cite{Pirce}, $\S \, 12.6.\quad \square$

\smallskip

Noether-Skolem's theorem  also implies that $\Gr_{k,\, l}'$ is
a homogeneous space of the group $\PGL_{kl}(\mathbb{C})$
represented as follows:
\begin{equation}
\nonumber
\Gr_{k,\, l}'=\PGL_{kl}(\mathbb{C})\bigl/
\PGL_k(\mathbb{C})\otimes \PGL_l(\mathbb{C}),
\end{equation}
where $\PGL_k(\mathbb{C})\otimes \PGL_l(\mathbb{C})$
denotes the image of the embedding
$\PGL_k(\mathbb{C})\times
\PGL_l(\mathbb{C})\hookrightarrow \PGL_{kl}(\mathbb{C})$
induced by the Kronecker product of matrices.

There is a {\it tautological} $M_k(\mathbb{C})$-bundle over the
matrix Grassmannian $\Gr'_{k,\, l}$ whose total space ${\cal A}_{k,\, l}'$
is defined in the following way:
$$
{\cal A}_{k,\, l}':=\{ (x,\, T)\mid x \in \Gr_{k,\, l}',\:
T\in M_{k,\, x}\}
\subset \Gr_{k,\, l}'\times M_{kl}(\mathbb{C}),
$$
where $M_{k,\, x}\subset M_{kl}(\mathbb{C})$ denotes the $k$-subalgebra corresponding to
a point $x\in \Gr_{k,\, l}'$, and the projection is induced by
the projection of the trivial bundle onto the first factor.

The noncompact space $\Gr_{k,\, l}'$ can be replaced
by a homotopy equivalent compact one
$\Gr_{k,\, l}.$ More precisely, consider the
``standard'' Hermitian metric $\langle A,\, B\rangle =
\tr(A\overline{B}^t)$ on $M_{kl}(\mathbb{C})$. We say that a $k$-subalgebra
$M_{k,\, x}\subset M_{kl}(\mathbb{C})$ is {\it unitary} if
$M_{k,\, x}=g(M_k(\mathbb{C}){\mathop{\otimes}
\limits_{\mathbb{C}}}\mathbb{C}E_l)g^{-1}$ for some $g\in \U(kl)$, where
$M_k(\mathbb{C}){\mathop{\otimes}
\limits_{\mathbb{C}}}\mathbb{C}E_l\subset M_{kl}(\mathbb{C})$ is the
``standard'' $k$-subalgebra.

Let $\PU(n)$ be the projective unitary group,
i.e. the quotient group
$\U(n)/\{ \exp(i\varphi)E_n
\mid \varphi \in \mathbb{R}\}.$ Put
\begin{equation}
\nonumber
\Gr_{k,\, l}:=\PU(kl)\bigl/
\PU(k)\otimes \PU(l).
\end{equation}
Obviously, it is a subspace in
$\Gr_{k,\, l}'$ which parameterizes unitary $k$-subalgebras.

In just the same way as in the
noncompact case, one can define the
tautological $M_k(\mathbb{C})$-bundle
${\cal A}_{k,\, l}\subset \Gr_{k,\, l}\times M_{kl}(\mathbb{C})$
over $\Gr_{k,\, l}$.

Recall that the group $\U(n)$ ($\PU(n)$) is a strong deformation
retract of $\GL_n(\mathbb{C})$ ($\PGL_n(\mathbb{C})$
respectively). Hence there is a homotopy equivalence $\Gr_{k,\,
l}\simeq \Gr_{k,\, l}'$. Because of this equivalence, we shall not
distinguish these spaces below.

\subsection{Floating algebra bundles}

The classical Grassmannians are classifying
spaces for vector bundles. In this subsection we introduce
a class of bundles that are classified by
{\it matrix} Grassmannians $\Gr_{k,\, l}$ (under the extra
condition $(k,\, l)=1$ whose sense will be clarified below).

Let $X$ be a finite $CW$-complex. By
$\widetilde M_n$ denote a trivial bundle
(over $X$) with fiber
$M_n(\mathbb{C})$ (note that in general a trivialization
on $\widetilde M_n$ is not supposed to be given).

\begin{definition}
\label{1}
Let $A_k$ $(k>1)$ be a locally trivial bundle over $X$
with fiber $M_k(\mathbb{C}).$
Assume that
there is a bundle map $\mu$
 \begin{equation}
 \nonumber
 \begin{array}{c}
 \diagram
 A_k\rrto^\mu \drto && \widetilde M_{kl}\dlto \\
 &X&
 \enddiagram
 \end{array}
 \end{equation}
such that for any point $x\in X$ it embeds the fiber
$(A_k)_x\cong
M_k(\mathbb{C})$ into the fiber $({\widetilde M}_{kl})_x\cong
M_{kl}(\mathbb{C})$
as a central simple subalgebra,
and the positive integers $k,\, l$ are relatively
prime (i.e. their greatest common
divisor $(k,l)~=~1$~)~.
Then the triple $(A_k,\, \mu ,\, \widetilde M_{kl})$ is called
a {\it floating algebra bundle} (abbrev. FAB) over $X$.
The locally trivial bundle $A_k$ is called a {\it core} of the FAB
$(A_k,\, \mu ,\, \widetilde M_{kl})$.
\end{definition}

\begin{remark}
Let $A$ be a central simple algebra over a field $\mathbb{K},$
$B\subset A$ a central simple subalgebra in $A$.
It is well known that the centralizer $\Z_A(B)$
of $B$ in $A$ is a central simple subalgebra in $A$ again,
moreover, the equality $A=B
{\mathop{\otimes}\limits_{\mathbb{K}}}\Z_A(B)$ holds \cite{Pirce}, $\S \, 12.7$.
Taking centralizers for all fibers
of the subbundle $A_k\subset \widetilde M_{kl}$
in the corresponding fibers of
the trivial bundle $\widetilde M_{kl}$,
we get the complementary subbundle $B_l$ with fiber $M_l(\mathbb{C})$
together with its embedding $\nu \colon
B_l\hookrightarrow \widetilde{M}_{kl}.$
Moreover, $A_k\otimes B_l=\widetilde M_{kl}.$

Conversely, to a given pair
$(A_k,\: B_l)$ consisting of
$M_k(\mathbb{C})$-bundle $A_k$ and $M_l(\mathbb{C})$-bundle
$B_l$ over $X$
such that $A_k\otimes B_l=\widetilde M_{kl},$ we can
construct a unique triple $(A_k,\, \mu ,\, \widetilde M_{kl}),$
where $\mu$ is the embedding
$A_k\hookrightarrow A_k\otimes B_l,\; a\mapsto a\otimes 1_{B_l}$.
\end{remark}

\begin{definition}
A {\it morphism from a FAB}
$(A_k,\, \mu ,\, \widetilde M_{kl})$ {\it to a FAB}
$(C_m,\, \nu ,\, \widetilde M_{mn})$
over $X$ is a pair $(f,\, g)$ of bundle maps
$f\colon A_k\hookrightarrow C_m,\;
g\colon \widetilde M_{kl}\hookrightarrow \widetilde M_{mn}$
such that
\begin{itemize}
\item $f,\, g$ are fiberwise homomorphisms of central algebras
(i.e. they actually are embeddings);
\item the square diagram
\begin{equation}
\nonumber
\begin{array}{c}
\diagram \widetilde M_{kl} \rto ^g & \widetilde M_{mn}
\\ A_k\uto ^\mu \rto ^f& C_m\uto _\nu
\enddiagram
\end{array}
\end{equation}
commutes;
\item let $B_l\subset \widetilde M_{kl},\;
D_n\subset \widetilde M_{mn}$ be the complementary
subbundles for $A_k,\; C_m$, respectively
(see the remark above), then
$g$ maps $B_l$ into $D_n$.
\end{itemize}
Note that a morphism $(f,\, g)\colon
(A_k,\, \mu ,\, \widetilde M_{kl})\rightarrow
(C_m,\, \nu ,\, \widetilde M_{mn})$ exists only if $k|m,\: l|n$.
\end{definition}

In particular,
an {\it isomorphism between FABs}
$(A_k,\, \mu ,\, \widetilde M_{kl})$ {\it and}
$(C_k,\, \nu ,\, \widetilde M_{kl})$ is a pair of bundle maps
$f:A_k\rightarrow C_k$, $g:\widetilde M_{kl}\rightarrow \widetilde
M_{kl}$ which are fiberwise isomorphisms of algebras
such that the diagram
\begin{equation}
\nonumber
\begin{array}{c}
\diagram \widetilde M_{kl} \rto ^g & \widetilde M_{kl}
\\ A_k\uto ^\mu \rto ^f& C_k\uto _\nu
\enddiagram
\end{array}
\end{equation}
commutes.

Clearly, FABs over $X$ with just defined
morphisms form a category $\mathfrak{FAB}(X)$.

For a continuous map $\varphi \colon X\rightarrow Y$
we have the natural transformation
$\varphi^* \colon \mathfrak{FAB}(Y)\rightarrow \mathfrak{FAB}(X)$.

Suppose $(k,l)=1.$ Put $\widetilde{\cal M}_{kl}:=\Gr_{k,\, l}
\times M_{kl}(\mathbb{C})$. Note that there is the {\it tautological FAB}
$({\cal A}_{k,\, l},\, \mu ,\, \widetilde
{\cal M}_{kl})$ over $\Gr_{k,\, l}$, where $\mu$ is the natural
inclusion existing by the definition of ${\cal A}_{k,\, l}$.

\begin{remark}
Note once more that there are compact and noncompact
versions of the considered theory. In the first case we consider
only unitary $k$-subalgebras and for any $x\in \Gr_{k,\, l}$
the embedding $\mu |_x \colon
({\cal A}_{k,\, l})_x\rightarrow M_{kl}(\mathbb{C})$ is a
unitary map (with respect to the standard Hermitian metrics).
\end{remark}

Let $\Psi_{k,\, l}(X)$ be the set of isomorphism
classes of FABs of the form
$(A_k,\, \mu,\, \widetilde{M}_{kl})$ over $X$. It is clear that
for such a FAB there is a classifying map $X\rightarrow \Gr_{k,\, l}$
(it can be obtained by fixing a trivialization on $\widetilde{M}_{kl}$).
Furthermore, the following assertion holds.

\begin{proposition}
\label{flocla}
{\rm (\cite{Prep}, Proposition 1.19)} If $\dim X<2\min \{ k,\, l\}$, then the assignment
$$
[X,\: \Gr_{k,\, l}]\rightarrow \Psi_{k,\, l}(X),
\quad \varphi \mapsto
\varphi^*({\cal A}_{k,\, l},\, \mu ,\, \widetilde
{\cal M}_{kl})
$$
is a bijection.
\end{proposition}

\subsection{The homogeneous space
$\Hom(M_k(\mathbb{C}),\, M_{kl}(\mathbb{C}))$}
\label{spaceofframes}

Fix a pair of coprime positive integers $k,\, l.$
Consider the set $\Hom(M_k(\mathbb{C}),\,
M_{kl}(\mathbb{C}))$ of all homomorphisms of
matrix algebras that take the unit $1_{M_k(\mathbb{C})}=E_k$
to the unit $1_{M_{kl}(\mathbb{C})}=E_{kl}$.
Since a matrix algebra is a simple ring, we see that
any such a homomorphism is an embedding.

There are two natural group actions on the set $\Hom(M_k(\mathbb{C}),\,
M_{kl}(\mathbb{C}))$. First, the group $\PGL_k(\mathbb{C})=
\Aut(M_k(\mathbb{C}))$ acts on $\Hom(M_k(\mathbb{C}),\,
M_{kl}(\mathbb{C}))$ as follows:
\begin{equation}
\label{action1}
\PGL_k(\mathbb{C})\times \Hom(M_k(\mathbb{C}),\,
M_{kl}(\mathbb{C})) \rightarrow \Hom(M_k(\mathbb{C}),\,
M_{kl}(\mathbb{C})),\; (g,\, \vartheta)\mapsto \vartheta \circ g^{-1}.
\end{equation}
Secondly, the group $\PGL_{kl}(\mathbb{C})=
\Aut(M_{kl}(\mathbb{C}))$ acts on $\Hom(M_k(\mathbb{C}),\,
M_{kl}(\mathbb{C}))$ in the following way:
\begin{equation}
\label{action2}
\PGL_{kl}(\mathbb{C})\times \Hom(M_k(\mathbb{C}),\,
M_{kl}(\mathbb{C})) \rightarrow \Hom(M_k(\mathbb{C}),\,
M_{kl}(\mathbb{C})),\; (g,\, \vartheta)\mapsto g\circ \vartheta.
\end{equation}
In particular, we have the corresponding action of the subgroup
$\PGL_{k}(\mathbb{C})=\PGL_{k}(\mathbb{C})\otimes E_l\subset
\PGL_{kl}(\mathbb{C}),$ where by $\otimes$ and $E_l$ we denote the Kronecker
product of matrices and the unit $l\times l$-matrix respectively.

Note that natural actions (\ref{action1}) and (\ref{action2})
do commute. Indeed, it is a general fact: to given sets $S$ and $R$
and a mapping $f\colon S \rightarrow R$ one has $(\alpha \circ f)\circ \beta
=\alpha \circ (f\circ \beta)$ for all $\alpha \in \Sigma(R),\:
\beta \in \Sigma(S),$ where $\Sigma(T)$ is the group of all bijections
of the set $T$.

We also can consider the set $\Hom(M_k(\mathbb{C}),\,
M_{kl}(\mathbb{C}))$ as a $\PGL_{k}(\mathbb{C})$ and
$\PGL_{kl}(\mathbb{C})$-invariant subset in the vector
space $\Hom_{Vect. sp}(L(M_{k}(\mathbb{C})),\, L(M_{kl}(\mathbb{C})))$,
where $L$ is the forgetful functor which to a matrix algebra
assigns its underlying vector space.
If we choose bases in $L(M_{k}(\mathbb{C}))$ and $L(M_{kl}(\mathbb{C}))$,
we can identify the latter space with the space of
$(kl)^2\times k^2$-matrices.
Then actions (\ref{action1}) and (\ref{action2}) correspond to
the right and left multiplication by invertible $k^2\times k^2$
and $(kl)^2\times (kl)^2$-matrices respectively.

Now let us describe $\Hom(M_k(\mathbb{C}),\,
M_{kl}(\mathbb{C}))$ as a homogeneous space
of the group $\PGL_{kl}(\mathbb{C})$.

\begin{definition}
A $k$-frame in $M_{kl}(\mathbb{C})$ is an ordered collection of
$k^2$ linearly independent matrices
$\{ \alpha_{i,\, j} \mid 1\leq i,\, j \leq k\}$
such that $\alpha_{i,\, j}
\alpha_{m,\, n}=\delta_{jm}\alpha_{i,\, n}$
for all $1\leq i,\, j,\, m,\, n\leq k,$ where $\delta_{jm}$
is the Kronecker delta-symbol, and
$\sum_{1\leq i\leq k}\alpha_{i,\, i}=E_{kl}$, where $E_{kl}\in
M_{kl}(\mathbb{C})$ is the unit matrix.
\end{definition}

Clearly, every such a $k$-frame is a basis in a certain (uniquely
determined by the frame) central $k$-subalgebra in
$M_{kl}(\mathbb{C}).$ For example, there is the ``standard''
$k$-frame $\{ e_{i,\, j}\mid 1\leq i,\, j\leq k\},$ where $e_{i,\,
j}:=E_{ij}\otimes E_l$ is the Kronecker product of a ``matrix
unit'' $E_{ij}$ of order $k$ by the unit $l\times l$ matrix $E_l$.
It is a frame in the subalgebra $M_k(\mathbb{C})
{\mathop{\otimes}\limits_{\mathbb{C}}} \mathbb{C}E_l\subset
M_{kl}(\mathbb{C}).$ Applying Noether-Skolem's theorem, we see that
all $k$-frames in $M_{kl}(\mathbb{C})$ are conjugate to each
other.

Thus, all $k$-frames in $M_{kl}(\mathbb{C})$
form the homogeneous space
$$
\Fr_{k,\, l}:=\PGL_{kl}(\mathbb{C})/E_k\otimes \PGL_l(\mathbb{C}).
$$
The notation means that the group $\PGL_{l}(\mathbb{C})$
is embedded into the group $\PGL_{kl}(\mathbb{C})$ by means of
the Kronecker product of matrices
$X\mapsto E_k\otimes X.$

\begin{proposition}
The space $\Hom(M_k(\mathbb{C}),\,
M_{kl}(\mathbb{C}))$ is a homogeneous space over the group
$\PGL_{kl}(\mathbb{C})$ with respect to the action (\ref{action2}).
Furthermore, there is an isomorphism of homogeneous spaces
$\Hom(M_k(\mathbb{C}),\, M_{kl}(\mathbb{C}))\cong \Fr_{k,\, l}$.
\end{proposition}
{\noindent \it Proof.\;} Fix a $k$-frame $\alpha$ in $M_{k}(\mathbb{C}),$
i.e. an ordered collection of
linearly independent matrices $\{ \alpha_{i,\, j}\mid 1\leq i,\, j
\leq k\},\; \alpha_{i,\, j}\in M_k(\mathbb{C})$
such that $\alpha_{i,\, j}\alpha_{r,\, s}=\delta_{jr}
\alpha_{i,\, s}.$ Then we have a one-to-one correspondence
between homomorphisms
$h\in \Hom(M_k(\mathbb{C}),\, M_{kl}(\mathbb{C}))$ and $k$-frames $\beta$
in $M_{kl}(\mathbb{C})$ given by $h\leftrightarrow \beta=
h(\alpha).\quad \square$

\smallskip

In particular, the previous proposition equips the set
$\Hom(M_k(\mathbb{C}),\,
M_{kl}(\mathbb{C}))$ with the topology of homogeneous space.

The space $\Fr_{k,\, l}$ is the total space of the principal
$\PGL_k(\mathbb{C})$-bundle
$$
\Fr_{k,\, l}\stackrel{\PGL_{k}(\mathbb{C})}{\longrightarrow}\Gr_{k,\, l}.
$$
The fiber of this bundle over $x\in \Gr_{k,\, l}$ consists of
all those $k$-frames that are contained in the $k$-subalgebra
$M_{k,\, x}\subset M_{kl}(\mathbb{C})$ corresponding to $x$.
The tautological $M_k(\mathbb{C})$-bundle
${\cal A}_{k,\, l}$ over $\Gr_{k,\, l}$ is associated with this
principal bundle.
Thus, the homogeneous space $\Fr_{k,\, l}$ can be treated as the total space
of the principal $\PGL_k(\mathbb{C})$-bundle
${\rm Prin}({\cal A}_{k,\, l})\rightarrow \Gr_{k,\, l}.$

\begin{remark}
There is also a unitary analog of this notion.
Clearly, the set of all {\it unitary} $k$-frames in $M_{kl}(\mathbb{C})$
(with respect to the standard Hermitian metric $\langle A,\, B\rangle =
\tr(A\overline{B}^t)$ on $M_{kl}(\mathbb{C})$)
is the homogeneous space
$\Fr_{k,\, l}:=\PU(kl)/E_k\otimes \PU(l).$ It is the
total space of the following
principal bundle
$$
\Fr_{k,\, l}\stackrel{\PU(k)}{\longrightarrow}\Gr_{k,\, l}
$$
over the (compact) matrix Grassmannian $\Gr_{k,\, l}.$
The fiber of this bundle over $x\in \Gr_{k,\, l}$
consists of all those {\it unitary} $k$-frames that are
contained in the $k$-subalgebra
$M_{k,\, x}\subset M_{kl}(\mathbb{C})$ corresponding to $x$
(recall that points $x\in \Gr_{k,\, l}$ correspond
to unitary central $k$-subalgebras in $M_{kl}(\mathbb{C})$).
The tautological $M_k(\mathbb{C})$-bundle
${\cal A}_{k,\, l}$ over $\Gr_{k,\, l}$ is associated with this principal bundle.
\end{remark}

Now note that there is a canonical FAB over $\Hom(M_k(\mathbb{C}),\,
M_{kl}(\mathbb{C}))=\Fr_{k,\, l}$ whose core (which is a trivial
$M_k(\mathbb{C})$-bundle equipped with the canonical trivialization)
is defined as the subspace
$$
\{ (h,\, h(T))\mid h\in \Hom(M_k(\mathbb{C}),\,
M_{kl}(\mathbb{C})),\: T\in M_k(\mathbb{C})\} \subset \Hom(M_k(\mathbb{C}),\,
M_{kl}(\mathbb{C}))\times M_{kl}(\mathbb{C}).
$$
Obviously, this FAB is the pull-back of the tautological FAB
$({\cal A}_{k,\, l},\, \mu,\, \widetilde{\cal M}_{kl})$ under the
projection $\Fr_{k,\, l}\rightarrow \Gr_{k,\, l}.$ Note that the
considered FAB over $\Hom(M_k(\mathbb{C}),\, M_{kl}(\mathbb{C}))$
is not trivial (see Definition \ref{triv2} below) although its
core is. The reason of this is that the embedding $\mu$ is ``twisted''.

Obviously, the assignment
$M_k(\mathbb{C})\mapsto \Hom(M_k(\mathbb{C}),\,
M_{kl}(\mathbb{C}))$ is functorial and therefore it can be
transferred to locally trivial bundles with fiber $M_k(\mathbb{C})$.
Thus, to any locally trivial $M_k(\mathbb{C})$-bundle $A_k$ over $X$
we assign the locally trivial $\Hom(M_k(\mathbb{C}),\,
M_{kl}(\mathbb{C}))$-bundle (over the same base $X$).
The obtained bundle we denote by ${\bf H}_{k,\, l}(A_k)$.
In particular, ${\bf H}_{k,\, l}(M_{k}(\mathbb{C}))\cong \Fr_{k,\, l},$
where $M_k(\mathbb{C})$ is regarded as a trivial bundle
over a point.

In terms of structure groups, the bundle ${\bf H}_{k,\, l}(A_k)$
is associated with the same principal $\PGL_k(\mathbb{C})$-bundle
as $A_k$, using the action (\ref{action1}).

\subsection{The $\Hom(M_k(\mathbb{C}),\,
M_{kl}(\mathbb{C}))$-bundle ${\bf H}_{k,\, l}(A_k^{univ})
\stackrel{p}{\rightarrow}\BPU(k)$}

By $A_k^{univ}$
denote the universal $M_k(\mathbb{C})$-bundle over $\BPU(k)$
(it is a locally trivial bundle with the structure group $\PU(k)
\simeq \PGL_k(\mathbb{C})=\Aut(M_k(\mathbb{C}))$).
In this subsection we study the $\Hom(M_k(\mathbb{C}),\,
M_{kl}(\mathbb{C}))$-bundle
${\bf H}_{k,\, l}(A_k^{univ})\stackrel{p}{\rightarrow}\BPU(k)$.

Let $t_{k,\, l}\colon \Gr_{k,\, l}\rightarrow \BPU(k)$
be the classifying map for ${\cal A}_{k,\, l}$
as a bundle with the structure group
$\PGL_k(\mathbb{C})$. Recall that the mapping
$t_{k,\, l}$ can be considered as a fibration in Homotopy Category.

First note that there is
the canonical embedding $p^*(A_k^{univ})\stackrel{\widetilde{\mu}}
{\hookrightarrow}
{\bf H}_{k,\, l}(A_k^{univ})\times M_{kl}(\mathbb{C})$
given by $(h,\, a)\mapsto (h,\, h(a)),$ where $a\in (A_k^{univ})_x,\:
h\in \Hom((A_k^{univ})_x,\, M_{kl}(\mathbb{C})).$ Thus the canonical
FAB $(p^*(A_k^{univ}),\, \widetilde{\mu},\,
{\bf H}_{k,\, l}(A_k^{univ})\times M_{kl}(\mathbb{C}))$ over
${\bf H}_{k,\, l}(A_k^{univ})$ is defined.

\begin{theorem}
\label{flid}
The total space ${\bf H}_{k,\, l}(A_k^{univ})$ is homotopy equivalent
to the matrix Grassmannian $\Gr_{k,\, l}.$ Moreover, the
bundle projection ${\bf H}_{k,\, l}(A_k^{univ})
\stackrel{p}{\rightarrow}\BPU(k)$
can be identified with the classifying map $t_{k,\, l}$
under this homotopy equivalence.
\end{theorem}

This theorem is a consequence of the following proposition.

Let $X$ be a finite $CW$-complex.
Consider the set $S_{k,\, l}(X)$ of equivalence classes of FABs of the form
$(A_k,\, \mu,\, \widetilde{M}_{kl})$ over $X$ with respect to the
following equivalence relation:
$$
(A_k,\, \mu,\, \widetilde{M}_{kl})\sim
(A'_k,\, \mu',\, \widetilde{M}_{kl}) \Leftrightarrow
$$
$$
\hbox{ there is a fiberwise homotopy } M \colon A_k\times I\rightarrow
\widetilde{M}_{kl} \hbox{ such that }
$$
$$
M |_{A_k\times \{ 0\} }=\mu,\; \, M |_{A_k\times \{ 1\} }\colon
A_k\stackrel{\cong}{\rightarrow}\mu'(A_k')\subset \widetilde{M}_{kl}
$$
(in particular, $A_k\cong A_k'$ and moreover
$(A_k,\, \mu,\, \widetilde{M}_{kl})\cong
(A'_k,\, \mu',\, \widetilde{M}_{kl})$).

\begin{proposition}
The spaces $\Gr_{k,\, l}$ and
${\bf H}_{k,\, l}(A_k^{univ})$ both
represent the homotopy functor $X\mapsto S_{k,\, l}(X)$.
\end{proposition}
{\noindent \it Proof.\;}
For $\Gr_{k,\, l}$ this is clear. Consider ${\bf H}_{k,\, l}(A_k^{univ}).$
Suppose $(A_k,\, \mu,\, \widetilde{M}_{kl})
\sim (A_k',\, \mu',\, \widetilde{M}_{kl})$
and $\varphi_0,\, \varphi_1\colon X\rightarrow
{\bf H}_{k,\, l}(A_k^{univ})$ are classifying maps for
$(A_k,\, \mu,\, \widetilde{M}_{kl}),\;
(A'_k,\, \mu',\, \widetilde{M}_{kl})$ respectively.
We have to construct a
homotopy $\Phi$ between $\varphi_0$ and $\varphi_1.$

Put $\bar{\varphi}_i:=p\circ \varphi_i,\, i=0,\, 1,$ where $p\colon
{\bf H}_{k,\, l}(A_k^{univ})\rightarrow \BPU(k)$ is the projection.
Since $A_k\cong A_k'$, we have a homotopy $\bar{\Phi}\colon
X\times I\rightarrow \BPU(k)$ between $\bar{\varphi}_0$ and
$\bar{\varphi}_1$. This (together with the covering homotopy property)
shows that without loss of generality we
can assume that $\bar{\varphi}_0=\bar{\varphi}_1=:\bar{\varphi}.$
In other words,
$A_k=A_k'=\bar{\varphi}^*(A_k^{univ})$ and $M$ defines a homotopy between
$\mu \colon A_k\rightarrow \widetilde{M}_{kl}$ and
$\mu'\colon A_k\rightarrow \widetilde{M}_{kl}$
(indeed, such a homotopy can be chosen as the composition
$M\circ \alpha,$ where $\alpha \colon A_k\stackrel{\cong}{\rightarrow}A_k'$
is an arbitrary isomorphism).
We have the bundle map $\widetilde{\varphi}\colon
{\bf H}_{k,\, l}(A_k)\rightarrow
{\bf H}_{k,\, l}(A^{univ}_k)$ which covers
$\bar{\varphi}\colon X\rightarrow \BPU(k).$
There is the natural one-to-one correspondence between sections
of the bundle ${\bf H}_{k,\, l}(A_k)\rightarrow X$ and embeddings
$A_k\hookrightarrow \widetilde{M}_{kl}=X\times M_{kl}(\mathbb{C})$
and also between corresponding homotopies.
Let $\sigma_i,\: i=0,\, 1$ be the sections corresponding to $\mu,\, \mu'$
respectively.
Therefore the homotopy $M$ gives us a homotopy $\varsigma$
between $\sigma_0$ and $\sigma_1$. Finally, the composite map
$\widetilde{\varphi}\circ \varsigma \colon X\times I\rightarrow
{\bf H}_{k,\, l}(A_k^{univ})$ is the required homotopy
between the classifying maps $\varphi_0\colon X\rightarrow
{\bf H}_{k,\, l}(A_k^{univ})$ and $\varphi_1\colon X\rightarrow
{\bf H}_{k,\, l}(A_k^{univ})$.

Now the converse assertion is clear.

Since the both spaces represent the same homotopy functor, there is
a homotopy equivalence
${\bf H}_{k,\, l}(A_k^{univ})\simeq \Gr_{k,\, l}.\quad \square$

\smallskip

Thus, we see that a lifting in the bundle ${\bf H}_{k,\, l}(A_k^{univ})\rightarrow \BPU(k)$
and an embedding of a locally trivial algebra bundle $A_k$ into a trivial one
(with fiber $M_{kl}(\mathbb{C})$) are the same things.
This will be useful for the study of obstructions for such an embedding in Subsection 3.2.

Now let us describe the homotopy equivalence
$\Gr_{k,\, l}\simeq {\bf H}_{k,\, l}(A_k^{univ})$
more explicitly. Obviously, the canonical embedding
$$
p^*(A_k^{univ})\stackrel{\widetilde{\mu}}
{\hookrightarrow}
{\bf H}_{k,\, l}(A_k^{univ})\times M_{kl}(\mathbb{C})
$$
determines the map
$\varphi \colon {\bf H}_{k,\, l}(A_k^{univ})\rightarrow \Gr_{k,\, l}$
such that $\varphi^*({\cal A}_{k,\, l},\, \mu,\, \Gr_{k,\, l}
\times M_{kl}(\mathbb{C}))=(p^*(A_k^{univ}),\, \widetilde{\mu},\,
{\bf H}_{k,\, l}(A_k^{univ})\times M_{kl}(\mathbb{C})),$ where
$({\cal A}_{k,\, l},\, \mu,\, \Gr_{k,\, l}
\times M_{kl}(\mathbb{C}))$ is the tautological FAB over $\Gr_{k,\, l}$.
(Indeed, $\varphi$ can be defined as follows. To a given
$h\in {\bf H}_{k,\, l}(A_k^{univ})$ the image
$\widetilde{\mu}(p^*(A_k^{univ})_h)\subset
M_{kl}(\mathbb{C})\cong \{ h\} \times M_{kl}(\mathbb{C})$ is a
subalgebra isomorphic to $M_{k}(\mathbb{C})$. We put $\varphi(h)=\alpha \in
\Gr_{k,\, l},$ where $\alpha$ is the point corresponding to the
$k$-subalgebra $\widetilde{\mu}(p^*(A_k^{univ})_h)\subset
M_{kl}(\mathbb{C})$.)
Since $t_{k,\, l}\colon \Gr_{k,\, l}\rightarrow \BPU(k)$
is a classifying map for ${\cal A}_{k,\, l}$, we have
$p\simeq t_{k,\, l}\circ \varphi.$

From the other hand, we have the fibration
$$
\Fr_{k,\, l}\rightarrow \Gr_{k,\, l}
\stackrel{t_{k,\, l}}{\longrightarrow}\BPU(k).
$$
The embedding $\mu \colon {\cal A}_{k,\, l}
\hookrightarrow \Gr_{k,\, l}\times M_{kl}(\mathbb{C})$
gives us a lifting $\psi \colon \Gr_{k,\, l}\rightarrow
{\bf H}_{k,\, l}(A_k^{univ})$ of $t_{k,\, l}.$
(Indeed, for $\alpha \in \Gr_{k,\, l},\: t_{k,\, l}(\alpha)=x$
we put $\psi(\alpha)=h,$ where $h\in ({\bf H}_{k,\, l}(A_k^{univ}))_x=
\Hom((A_k^{univ})_x,\, M_{kl}(\mathbb{C}))=\Hom(({\cal A}_{k,\, l})_{\alpha},
\, M_{kl}(\mathbb{C}))$ is the homomorphism defined by $\mu |_{\alpha}$).
Clearly,
$\psi^*(p^*(A_k^{univ}),\, \widetilde{\mu},\,
{\bf H}_{k,\, l}(A_k^{univ})\times M_{kl}(\mathbb{C}))=
({\cal A}_{k,\, l},\, \mu,\, \widetilde{\cal M}_{kl})$. In particular,
$p\circ \psi \simeq t_{k,\, l}.$

We have already observed that
since the both spaces ${\bf H}_{k,\, l}(A_k^{univ})$ and $\Gr_{k,\, l}$
represent the same homotopy functor, they are homotopy equivalent
to each other.
Moreover it is easy to see that the above constructed
maps $\varphi \colon {\bf H}_{k,\, l}(A_k^{univ})\rightarrow \Gr_{k,\, l}$ and
$\psi \colon \Gr_{k,\, l}\rightarrow
{\bf H}_{k,\, l}(A_k^{univ})$ are the desired
homotopy equivalences,
$\psi \circ \varphi \simeq \id_{{\bf H}_{k,\, l}(A_k^{univ})},\;
\varphi \circ \psi \simeq \id_{\Gr_{k,\, l}}.\quad \square$

\smallskip

Thus, we can identify the fibrations (in Homotopy Category)
$$
\Fr_{k,\, l}\rightarrow \Gr_{k,\, l}
\stackrel{t_{k,\, l}}{\longrightarrow}\BPU(k)
$$
and
$$
\Fr_{k,\, l}\rightarrow {\bf H}_{k,\, l}(A_k^{univ})
\stackrel{p}{\rightarrow}\BPU(k).
$$

\section{The functor represented by the $H$-space $\Gr$}

In this section we give a brief survey of the stable
theory of floating algebra bundles. We do this because
this theory is needed to motivate the definition of the generalized
Brauer group in Subsection 3.3 (see for example Theorem \ref{thimportfib} and its corollary).
For more details see \cite{Prep}, Ch.2.

\subsection{The stable equivalence of FABs}

Define the product $\circ$ of two
FABs over $X$
$(A_k,\, \mu ,\, \widetilde{M}_{kl}) ,\;
(B_m,\, \nu ,\, \widetilde{M}_{mn})$
such that $(km,\, ln)=1$ as
$$
(A_k,\, \mu ,\, \widetilde{M}_{kl})
\circ (B_m,\, \nu ,\, \widetilde{M}_{mn}):=
(A_k\otimes B_m,\, \mu \otimes \nu,\, \widetilde{M}_{kl}
\otimes \widetilde{M}_{mn})
$$
(notice that $\widetilde{M}_{kl}
\otimes \widetilde{M}_{mn}=\widetilde{M}_{klmn}$).
\begin{definition}
\label{triv2}
A FAB of the form
$(\widetilde M_k,\, \tau ,\, \widetilde M_{kl})$ is called
{\it trivial} if
the map $\tau \colon \widetilde{M}_k
\rightarrow \widetilde{M}_{kl}$ is the following:
$$
X\times
M_k(\mathbb{C})\rightarrow X\times M_{kl}(\mathbb{C}),\qquad
(x,\, T)\mapsto (x,\, T\otimes E_l) \; \hbox{ for all } x\in X
$$
(under some choice of trivializations
on $\widetilde{M}_k$ and $\widetilde{M}_{kl}$),
where $E_l$ is the unit $l\times
l$-matrix and $T\otimes E_l$ denotes the Kronecker product
of matrices.
In other words, the bundle $\widetilde M_k$ is embedded
into $\widetilde M_{kl}$ as a fixed subalgebra.
\end{definition}

\begin{definition}
\label{floating}
Two FABs $(A_k,\, \mu ,\, \widetilde M_{kl})$ and
$(B_m,\, \nu ,\, \widetilde
M_{mn})$ over $X$ are said to be {\it stable equivalent}
if there is a sequence of pairs
$\{t_i,\, u_i\}\in \mathbb{N}^2$,
$1\leq i\leq s$ such that
\begin{itemize}
\item \; $\{t_1,\, u_1\}=\{k,\, l\},\; \{t_s,\, u_s\}=\{m,\, n\};$
\item \; $(t_it_{i+1},\, u_iu_{i+1})=1$
if $s>1,\, 1\leq i\leq s-1,$
\end{itemize}
and a corresponding sequence of FABs $(A_{t_i},\, \mu _i,\,
\widetilde M_{t_iu_i})$ over $X$ such that
\begin{itemize}
\item \; $(A_{t_1},\, \mu _1,\, \widetilde M_{t_1u_1})=
(A_k,\, \mu ,\, \widetilde
M_{kl}),\; (A_{t_s},\, \mu _s,\, \widetilde M_{t_su_s})=
(B_m,\, \nu, \,
\widetilde M_{mn});$
\item \; $(A_{t_i},\, \mu _i,\,
\widetilde M_{t_iu_i})\circ (\widetilde
M_{t_{i+1}},\, \tau ,\, \widetilde M_{t_{i+1}u_{i+1}})
\cong (A_{t_{i+1}}
,\, \mu_{i+1},\, \widetilde M_{t_{i+1}u_{i+1}})\circ
(\widetilde M_{t_i},\, \tau,\, \widetilde M_{t_iu_i})$,
\end{itemize}
where $1\leq i\leq s-1$ and
$(\widetilde M_{t_i},\, \tau ,\, \widetilde
M_{t_iu_i})$ are trivial FABs.
\end{definition}

By $\widetilde{\AB}^1(X)$ denote the set of stable equivalence
classes of FABs over $X$.

The following theorem justifies
the previous definition. Note that a homomorphism
of central algebras $M_{kl}(\mathbb{C})\hookrightarrow
M_{klmn}(\mathbb{C})$ induces the corresponding map of matrix
Grassmannians $i_{k,\,l;\, m,\, n}\colon \Gr_{k,\, l}
\rightarrow \Gr_{km,\, ln}$; the direct limits below are taken over such maps.

\begin{theorem}
\label{homeq}
{\rm ( \cite{Prep}, Proposition 2.1)} 1) For all sequences of pairs of positive
integers $\{ k_j,\, l_j\}_{j\in\mathbb{N}}$
such that
 $$
 {\rm(i)} \quad k_j,\: l_j\to\infty;\quad{\rm(ii)}\ k_j|k_{j+1},\
 l_j|l_{j+1};\quad{\rm(iii)}\ (k_j,\, l_j)=1
 $$
for every $j,$ the corresponding direct limits
$\varinjlim_j \Gr_{k_j,\, l_j}$ are homotopy-equivalent.
This unique homotopy type
we denote by $\Gr$.\\
2) $\Gr$ is a classifying space for stable equivalence classes
of FABs over a finite $CW$-complex $X.$ In other words,
the functor $X\mapsto \widetilde{\AB}^1(X)$
from the homotopy category of finite $CW$-complexes
to the category $\mathfrak{Sets}$ is represented
by the space $\Gr$.
\end{theorem}
The proof is based on Proposition \ref{flocla}
and on the following lemma.

\begin{lemma}
If $(km,ln)=1$, then
the embedding
$$
i_{k,\,l;\, m,\, n}\colon \Gr_{k,\, l}\rightarrow \Gr_{km,\, ln}
$$
is a homotopy equivalence in dimensions $<2\min \{k,\, l\}$.
\end{lemma}
Therefore for any finite $CW$-complex $X$,
$\dim(X)<2\min \{k,\, l,\, m,\, n\}$ and for any
map $\varphi_{km,\, ln}\colon X\rightarrow \Gr_{km,\, ln}$
there are maps $\varphi_{k,\, l}\colon X\rightarrow \Gr_{k,\, l}$
and $\varphi_{m,\, n}\colon X\rightarrow \Gr_{m,\, n}$
such that $i_{k,\,l;\, m,\, n}\circ \varphi_{k,\, l}\simeq
\varphi_{km,\, ln}\simeq i_{m,\,n;\, k,\, l}\circ \varphi_{m,\, n}.$
Note also that $i_{k,\,l;\, m,\, n}^*
({\cal A}_{km,\, ln},\, \mu \, ,\widetilde{\cal M}_{klmn})=
({\cal A}_{k,\, l},\, \mu \, ,\widetilde {\cal M}_{kl})\circ
(\widetilde {\cal M}_{m},\, \tau,\, \widetilde {\cal M}_{mn})$,
whence the stable equivalence relation.

\subsection{The group structure}

Let $(A_k,\, \mu ,\, \widetilde{M}_{kl})$ be a FAB over $X$.
By $[(A_k,\, \mu ,\, \widetilde{M}_{kl})]$ we denote its
stable equivalence class (with respect to the equivalence relation defined
in the previous subsection).
Define the product $\diamond$ of two classes
$[(A_k,\, \mu ,\, \widetilde{M}_{kl})] ,\;
[(B_m,\, \nu ,\, \widetilde{M}_{mn})]$ for $(km,\, ln)=1$
as
$$
[(A_k,\, \mu ,\, \widetilde{M}_{kl})]
\diamond [(B_m,\, \nu ,\, \widetilde{M}_{mn})]=
[(A_k,\, \mu ,\, \widetilde{M}_{kl})\circ
(B_m,\, \nu ,\, \widetilde{M}_{mn})]
$$
$$=[(A_k\otimes B_m,\, \mu \otimes \nu ,\, \widetilde{M}_{klmn})].
$$
Clearly, this product is well defined.
The following lemma allows us to reject the restriction $(km,\, ln)=1$.
\begin{lemma}
\label{repres}
For any pair $\{ k,\, l\}$ such that $\rm{(i)}$ $(k,l)~=1,$
$\rm{(ii)}$ $2\min\{k,\, l\}\geq \dim X,$
any stable equivalence class of FABs over
$X$ has a representative of the form
$(A_k,\, \mu ,\, \widetilde M_{kl}).$
\end{lemma}

Clearly, the product $\diamond$ is associative, commutative,
and has identity element
$[(\widetilde M_k,\, \tau ,\, \widetilde M_{kl})]$,
where $(\widetilde M_k,\, \tau ,\, \widetilde M_{kl})$
is a trivial FAB. Moreover, for any
class $[(A_k,\, \mu ,\, \widetilde{M}_{kl})]$ there
exists the inverse element. In order to find it, let us
recall the following fact.
The centralizer
$\Z_P(Q)$ of a central simple subalgebra $Q$
in a central simple algebra $P$ (over some field
$\mathbb{K}$) is a central simple subalgebra again,
moreover, the equality $P=Q
{\mathop{\otimes}\limits_{\mathbb{K}}}\Z_P(Q)$ holds.
Therefore by taking centralizers for every fiber
of the subbundle $A_k$ in $\widetilde M_{kl},$
we obtain the complementary subbundle
$B_l$ with fiber $M_l(\mathbb{C})$
together with its embedding $\nu \colon
B_l\hookrightarrow \widetilde{M}_{kl}$
into the trivial bundle. Moreover,
$A_k\otimes B_l=\widetilde M_{kl}.$
It is not hard to prove that
$[(B_l,\, \nu ,\, \widetilde{M}_{kl})]$ is the
inverse element for
$[(A_k,\, \mu ,\, \widetilde{M}_{kl})]$.
Thus, the functor $X\mapsto \widetilde{\AB}^1(X)$
takes values in the category of Abelian groups
$\mathfrak{Ab}$.

\begin{proposition}
{\rm (\cite{Prep}, Theorem 2.5)} The space $\Gr$ can be equipped with an $H$-space structure
such that there is a natural equivalence
between $X\mapsto [X,\: \Gr]$ and
$X\mapsto \widetilde{\AB}^1(X)$ as functors to
the category $\mathfrak{Ab}.$
\end{proposition}

\subsection{Some properties of FAB's core}

Recall that a locally trivial
$\Aut(M_k(\mathbb{C}))\cong \PGL_k(\mathbb{C})$ (or $\PU(k)$)-bundle
$A_k$ is called a {\it core} of a FAB
$(A_k,\, \mu ,\, \widetilde M_{kl})$ (if such a FAB exists, of course).
We mention its properties because they will play an important role in the
definition of the generalized Brauer group in Subsection 3.3.

\begin{lemma}
If $A_k$ is the core
of some FAB $(A_k,\, \mu ,\, \widetilde M_{kl})$, then
its structure group
can be reduced from $\Aut M_k(\mathbb{C})\cong \PGL_k(\mathbb{C})$
to $\SL_k(\mathbb{C})$ (or equivalent from $\PU(k)$ to $\SU(k)$).
\end{lemma}
{\noindent \it Proof.}\; Indeed, since $\rho_{kl}=\rho_k\times
\rho_l$ for $(k,\, l)=1$, where $\rho_n$ is the group of $n$th
degree roots of unity, we have
$$
\Gr_{k,\, l}=\PU(kl)/\PU(k)\otimes
\PU(l)=\SU(kl)/\SU(k)\otimes \SU(l).\quad \square
$$

The following lemma describes a characteristic property of cores.

\begin{lemma}
\label{biglem}
(\cite{Prep}, Lemma 2.7)
Let $X$ be a finite $CW$-complex.
Suppose $\dim X\leq 2\min\{k,\, m\};$
then the following conditions are equivalent:
\begin{itemize}
\item
$A_k$ is the core of some FAB over $X$;
\item
for arbitrary~$m$ such that $2m\geq \dim X$
there is a bundle $B_m$
with fiber $M_m(\mathbb{C})$ such that
$A_k\otimes{\widetilde M}_m\cong B_m\otimes {\widetilde M}_k;$
\item
$A_k\otimes {\widetilde M}_m\cong
B_m\otimes {\widetilde M}_k$
for some locally trivial bundle $B_m$ with fiber
$M_m(\mathbb{C})$, where $(k,m)=1$.
\end{itemize}
Moreover, for any pair of bundles $A_k,\: B_m$
such that $(k,m)=1$ and
$A_k\otimes {\widetilde M}_m\cong
B_m\otimes {\widetilde M}_k,$
there exists a unique stable equivalence class
of FABs over $X$ which has (for sufficiently large
$n,\: (km,n)=1$) FABs of the forms $(A_k,\, \mu ,\,
\widetilde M_{kn}),\;
(B_m,\, \nu ,\, \widetilde M_{mn})$ as representatives
(for some embeddings $\mu,\: \nu$).
\end{lemma}

\subsection{Localization}

Let $X$ be a finite $CW$-complex, $k\geq 2$
a fixed integer.
The set of isomorphism classes of bundles of
the form $A_{k^m}$ (for arbitrary $m\in \mathbb{N}$) over
$X$ with fiber $M_{k^m}(\mathbb{C})$ is
a monoid with respect to the operation $\otimes$
(with the identity element $M_{k^0}(\mathbb{C})\cong\mathbb{C}$).

Let us consider the following equivalence relation
$$
A_{k^m}\sim B_{k^n}
\Longleftrightarrow
\exists r,s\in\mathbb{N} \mbox{\: such that\:}
A_{k^m}
\otimes\widetilde M_{k^r}
\cong B_{k^n}
\otimes\widetilde M_{k^s}
$$
($\Rightarrow m+r=n+s$).
The set of equivalence classes $[A_{k^m}]$
of such bundles is a group with respect to the
operation induced by $\otimes$. This group we denote by
$\widetilde{\AB}^k(X).$

Let us consider the direct limit
$\lim\limits_{\scriptstyle
\longrightarrow\atop \scriptstyle n}
\BPU(k^n)$
with respect to the maps induced by
$$
\begin{array}{ccc}
\PU(k^n) & \hookrightarrow &
\PU(k^{n+1}), \\
A & \mapsto & E_k\otimes A.
\end{array}
$$
Clearly, the functor $X\mapsto \widetilde{\AB}^k(X)$
is represented by $\lim\limits_{\scriptstyle
\longrightarrow\atop \scriptstyle n}
\BPU(k^n).$

According to Lemma \ref{repres}, for any stable equivalence class
of FABs over $X$ there is a representative of the form
$(A_{k^m},\, \mu ,\, \widetilde{M}_{(kl)^m}),\; (k,l)=1.$
Therefore for any $k$ we have the group homomorphism
$$\widetilde{\AB}^1(X)\rightarrow
\widetilde{\AB}^k(X),\quad [(A_{k^m},\, \mu ,\,
\widetilde{M}_{(kl)^m})]\mapsto
[A_{k^m}]$$
induced by the following map of the direct limits
$$
 \diagram
 &\\
 \Gr_{k^2,\, l^2}\rto^{t_{k^2,\, l^2}}
 \uto & \BPU(k^2)\uto \\
 \Gr_{k,\, l}\uto^{i_{k,\, l;\, k,\, l}} \rto^{t_{k,\, l}} & \BPU(k),\uto
 \enddiagram
$$
where $t_{k^2,\, l^2}$ and $t_{k,\, l}$
are classifying maps for the cores
${\cal A}_{k,\, l}$ and ${\cal A}_{k^2,\, l^2}$
as $\PU(k)$ and $\PU(k^2)$-bundles,
respectively.

The kernel of the homomorphism
$\widetilde{\AB}^1(X)\rightarrow
\widetilde{\AB}^k(X)$ is just
the $k$-torsion subgroup in
$\widetilde{\AB}^1(X)$.

Set $t_k:=\varinjlim_rt_{k^r,\, l^r}\colon \Gr \rightarrow
\varinjlim_r \BPU(k^r)=:\BPU(k^\infty).$ Then $t_k$ is the
composition of the localization map (\cite{Sullivan}, Ch.2)
$\widetilde{t}_k\colon \Gr \rightarrow \varinjlim_r\BSU(k^r)=
\BSU[\frac{1}{k}]$ and the natural map
$\varinjlim_r\BSU(k^r)\rightarrow \varinjlim_r\BPU(k^r)$ induced
by the group epimorphisms $\SU(k^r)\rightarrow \PU(k^r)$ with the
kernels $\{ \lambda E_{k^r}\mid \lambda^{k^r}=1\}\cong
\rho_{k^r}.$

\subsection{Relation between $\widetilde{\AB}^1$ and
$\widetilde{\KSU}$-theory}

Recall that $\BSU_\otimes$ is the space $\BSU$ with
the structure of $H$-space related to
the tensor product of virtual $\SU$-bundles
of virtual dimension $1$.

\begin{theorem}
\label{Kthcoincid}
{\rm (\cite{Prep}, Theorem 2.17)}
There is an $H$-space isomorphism
$\Gr \cong \BSU_\otimes.$
\end{theorem}

By $\widetilde{\KSU}(X)$ denote the reduced $\K$-functor
constructed by means of $\SU$-bundles over $X$.
Recall that $\widetilde{\KSU}(X)$ is a ring with the
multiplication induced by the tensor product of bundles.

The previous theorem claims that the group $\widetilde{\AB}^1(X)$
is isomorphic to the multiplicative group of the ring
$\widetilde{\KSU}(X),$ i.e. the group (because $X$ is a finite
$CW$-complex) of elements of $\widetilde{\KSU}(X)$ with respect to
the operation $\xi \ast \eta =\xi +\eta +\xi \eta$ $(\xi, \eta \in
\widetilde{\KSU}(X)$, i.e. $\xi,\, \eta$ are of virtual dimension
$0$).

This gives us a geometric description
of the $H$-structure on $\BSU_\otimes.$
For example, the construction of the inverse stable equivalence class
$[(B_m,\, \nu,\, \widetilde{M}_{mn})]$
for a given one $[(A_k,\, \mu,\, \widetilde{M}_{kl})]$
is closely connected with
taking centralizer for a subalgebra in a fixed matrix
algebra.

Although the theory of
floating algebra bundles leads to the theory which can also be described in the classical terms, the obtained
geometric approach might be applicable in other branches of Topology (for example,
in Cobordism theory, cf. \cite{e1}, \cite{e2}).

\subsection{A $\U$-version}

Consider the canonical map $\BU(k)\rightarrow
\BPU(k)$ induced by the group homomorphism
$\U(k)\rightarrow \PU(k).$
By $\widehat{\Gr}_{k,\, l}$ denote the total space
of the $\Fr_{k,\, l}$-fibration (recall that $\Fr_{k,\, l}$
denotes the space of (unitary) $k$-frames
in $M_{kl}(\mathbb{C})$, see Subsection \ref{spaceofframes}) induced by
the fibration
$\Gr_{k,\, l}\stackrel{\Fr_{k,\, l}}{\longrightarrow}
\BPU(k)$ and the map
$\BU(k)\rightarrow
\BPU(k)$ (as ever, the integers $k,\, l$ are
assumed to be coprime), i.e. the fiber product
\begin{equation}
\label{unitv}
\begin{array}{c}
\diagram
\widehat{\Gr}_{k,\, l} \dto \rto^{\widehat{t}_{k,\, l}} & \BU(k) \dto \\
\Gr_{k,\, l} \rto^{t_{k,\, l}\; \;} & \BPU(k). \\
\enddiagram
\end{array}
\end{equation}

It follows easily that there is a $\mathbb{C}P^\infty$-fibration
$\widehat{\Gr}_{k,\, l}\rightarrow{\Gr}_{k,\, l}.$

\begin{remark}
\label{unitver}
Let us give a description of $\widehat{\Gr}_{k,\, l}$ analogous to the one
for $\Gr_{k,\, l}$ in Subsection 1.4.
Let $\xi_k^{univ}\rightarrow \BU(k)$ be the universal vector bundle
with fiber $\mathbb{C}^k.$ Then clearly
${\bf H}_{k,\, l}(\End (\xi_k^{univ}))\simeq \widehat{\Gr}_{k,\, l},$
moreover under this identification the projection
${\bf H}_{k,\, l}(\End (\xi_k^{univ}))
\rightarrow \BU(k)$ coincides with the map
$\widehat{\Gr}_{k,\, l}\rightarrow \BU(k)$ we have just
defined by the commutative square.
\end{remark}

Consider the following morphism of $\U(k)$-fibrations:
\begin{equation}
\label{unitver1}
\diagram
& \U(k) \rto & \EU(k) \dto \\
\U(k) \urto^= \rto & \Fr_{k,\, l} \dto \urto &
\BU(k) \\
& \widehat{\Gr}_{k,\, l}, \urto^{\widehat{t}_{k,\, l}} \\
\enddiagram
\end{equation}
where $\widehat{t}_{k,\, l}$ is the classifying map
for the canonical $\U(k)$-bundle over $\widehat{\Gr}_{k,\, l}$ and
by $\EU(k)$ we denote the total space of the universal principal
$\U(k)$-bundle which is contractible. A simple computation with
homotopy sequences of the fibrations shows that $\widehat{t}_{k,\,
l*}\colon \pi_{2r}(\widehat{\Gr}_{k,\, l})\rightarrow
\pi_{2r}(\BU(k)),\, r\leq \min \{k,\, l\}$ is just the
monomorphism $\mathbb{Z}\rightarrow \mathbb{Z},\; 1\mapsto k\cdot
1$ (note that the odd-dimensional stable homotopy groups of both
spaces are equal to $0$). This implies (\cite{Sullivan}, Theorem
2.1) that the direct limit map $\widehat{t}_{k}\colon
\widehat{\Gr}\rightarrow \varinjlim_r\BU(k^r)=:\BU(k^\infty)$ is
just the localization away from $k$ (in the sense that $k$ is
invertible; in particular, $\varinjlim_r\BU(k^r)=\BU[\frac{1}{k}]$
is a $\mathbb{Z}[\frac{1}{k}]$-local space), where
$\widehat{\Gr}:=\varinjlim_{(k,\, l)=1}\widehat{\Gr}_{k,\, l}.$

The diagram (\ref{unitv}) gives us the diagram
\begin{equation}
\label{uncasedia}
\diagram
\widehat{\Gr} \dto_{\mathbb{C}P^\infty} \rto^{\widehat{t}_k\quad} & \BU(k^\infty) \dto^{\mathbb{C}P^\infty} \\
\Gr \rto^{t_{k}\quad \; \;} & \BPU(k^\infty). \\
\enddiagram
\end{equation}
The left-hand vertical arrow in terms of bundles can be described as an assignment
$$
(\xi_{k^n},\, (\End(\xi_{k^n}), \, \mu,\, \widetilde{M}_{(kl)^n}))\mapsto
(\End(\xi_{k^n}), \, \mu,\, \widetilde{M}_{(kl)^n}),
$$
$\widehat{t}_k$ as an assignment
$$
(\xi_{k^n},\, (\End(\xi_{k^n}), \, \mu,\, \widetilde{M}_{(kl)^n}))\mapsto \xi_{k^n},
$$
and also $t_k$ and the right-hand arrow as $(A_{k^n}, \, \mu,\, \widetilde{M}_{(kl)^n})\mapsto A_{k^n}$
and $\xi_{k^n}\mapsto \End(\xi_{k^n})$ respectively, where
$(\End(\xi_{k^n}), \, \mu,\, \widetilde{M}_{(kl)^n})$ and $(A_{k^n}, \, \mu,\, \widetilde{M}_{(kl)^n})$
denote some FABs.

The space $\widehat{\Gr}$ is an $H$-space
with respect to the multiplication
induced by the tensor product of bundles. It can be proved
that $\widehat{\Gr}\cong \BU_{\otimes}$ as $H$-spaces.
Let us also recall that $\BU_{\otimes}\cong \BSU_{\otimes}
\times \mathbb{C}P^{\infty}$ and $\Gr \cong \BSU_{\otimes}$
as $H$-spaces, hence $\widehat{\Gr}\cong \Gr \times
\mathbb{C}P^{\infty}.$ In particular, the $H$-space
$\widehat{\Gr}$ represents the functor of the
``multiplicative group'' of the ring
$\widetilde{\K}_{\mathbb{C}},$
i.e. the functor $X\mapsto \widetilde{\K}_{\mathbb{C}}(X),$
where $\widetilde{\K}_{\mathbb{C}}(X)$ is considered as a group
with respect to the operation $\xi * \eta=\xi +\eta +
\xi \eta,\; \xi,\, \eta \in \widetilde{\K}_{\mathbb{C}}(X)$
(here $\widetilde{\K}_{\mathbb{C}}$ is the reduced
complex $K$-functor).

By $\widehat{\AB}^1(X),\; \widehat{\AB}^k(X)$ denote the groups $[X,\, \widehat{\Gr}]$ and $[X,\, \BU(k^\infty)]$
respectively. We also have the group homomorphism $\widehat{t}_{k*}(X)\colon
\widehat{\AB}^1(X)\rightarrow \widehat{\AB}^k(X)$ induced by $\widehat{t}_k$.

\section{Homotopy invariants related to algebra bundles}

\subsection{Reminder: The classical Brauer group}

Undoubtedly, the most important invariant constructed
by means of matrix algebra bundles (called in this context ``Azumaya bundles'')
is the Brauer group which plays an
important role not only in Topology
but also in Algebraic Geometry (where it turns out to be
a birational invariant of varieties \cite{Grothendieck}). In this
subsection we give a brief survey of the classical results concerning
its ``purely topological'' version.

So let $X$ be a finite $CW$-complex.
Consider the set of isomorphism classes
of locally trivial bundles over $X$
with fiber $M_k(\mathbb{C})$ with an arbitrary
integer $k>1$. On the set of such bundles
consider the following equivalence relation:
\begin{equation}
\label{classicequiv}
A_k\sim B_l\; \Leftrightarrow A_k\otimes \widetilde{M}_m
\cong B_l\otimes \widetilde{M}_n \hbox{\quad for some }m,\, n>1,
\end{equation}
where by $\widetilde{M}_r$
we denote a trivial bundle over $X$ with fiber
$M_r(\mathbb{C})$. By $\{ C_m\}$ denote
the stable equivalence class of $C_m$.
It can easily be checked that the product
$\{ A_k\} \circ \{ B_l\} :=\{ A_k\otimes B_l\}$
is well defined and equips the set of
stable equivalence classes
with the structure of Abelian group.
We denote this group by $\widetilde{\AB}(X).$ Clearly,
$\widetilde{\AB}(X)=\varinjlim_k\widetilde{\AB}^k(X),$ where
the direct limit is taken over the natural homomorphisms
$\widetilde{\AB}^k(X)\rightarrow \widetilde{\AB}^{km}(X)$ induced by the assignments
$A_{k^n}\mapsto A_{(km)^n}:=A_{k^n}\otimes \widetilde{M}_{m^n}$.

Now consider the following more coarse stable equivalence relation:
$$
A_k\sim B_l\; \Leftrightarrow \exists \, \xi_m,\,
\eta_n \hbox{ such that }A_k\otimes \End(\xi_m)
\cong B_l\otimes \End(\eta_n),
$$
where by $\xi_m,\, \eta_n$ we denote vector bundles over $X$ of rank
$m,\, n$, respectively.
By $[C_m]$ denote the stable equivalence
class of $M_m(\mathbb{C})$-bundle $C_m.$
The tensor product of bundles
induces a group structure on the set of such stable equivalence
classes. It is just the classical topological
Brauer group $\Br(X).$

Let us give a homotopic description of $\Br(X).$
Consider the fibration
\begin{equation}
\label{obstruc}
\begin{array}{ccc}
\mathbb{C}P^{\infty} & \hookrightarrow & \BU(k) \\
&& \downarrow \\
&& \BPU(k) \\
\end{array}
\end{equation}
corresponding to the exact sequence of groups
$$
1\rightarrow \U(1) \rightarrow \U(k)
\rightarrow \PU(k)\rightarrow 1.
$$
The first obstruction $\alpha(f)$ for the lifting of a map
$f\colon X\rightarrow \BPU(k)$ in
(\ref{obstruc}) belongs to $H^3(X;\, \mathbb{Z})$ (\cite{Fuchs}, Ch.5)
and has order $k$ (the last assertion follows from
the fact that $\alpha(f)=\delta(\beta(f)),$ where
$\beta(f)\in H^2(X;\, \mathbb{Z}/k\mathbb{Z})$ is the first
obstruction for the lifting in the bundle $\B \mu_k \hookrightarrow
\BSU(k)\rightarrow \BPU(k)$
and $\delta \colon H^2(X;\, \mathbb{Z}/k\mathbb{Z})\rightarrow
H^3(X;\, \mathbb{Z})$ is the coboundary homomorphism).

Let us introduce the following notation for the direct limits:
$$
\B{\cal U}:=
\varinjlim_k\BU(k),\quad
\B{\cal PU}:=
\varinjlim_k\BPU(k),
$$
where $k$ runs over all positive integers and the limits are taken over
the maps inducing by the tensor product
(in particular, $\B{\cal U}$ is a $\mathbb{Q}$-space).
Recall that
$$
\B {\cal U}=\prod_{q\geq 1}
\K(2q,\, \mathbb{Q}),\; \;
\B {\cal PU}=\K(2,\, \mathbb{Q}/\mathbb{Z})\times
\prod_{q\geq 2} \K(2q,\, \mathbb{Q}),\; \;
\mathbb{C}P^\infty =\K(2,\, \mathbb{Z}),
$$
and we have the fibration
$$
\mathbb{C}P^\infty \hookrightarrow
\K(2,\, \mathbb{Q})\rightarrow \K(2,\, \mathbb{Q}/\mathbb{Z})
$$
corresponding to the exact sequence of the coefficients
groups
$$
0\rightarrow \mathbb{Z}\rightarrow \mathbb{Q} \rightarrow
\mathbb{Q}/\mathbb{Z}\rightarrow 0.
$$
We consider every space $\mathbb{C}P^\infty,\,
\B {\cal U},$ and $\B {\cal PU}$
as an $H$-space with the multiplication induced by
the tensor product of bundles which it classifies
(for example, $\mathbb{C}P^\infty$ classifies complex line bundles,
and there is the isomorphism $\mathbb{C}P^\infty \cong \K(2,\, \mathbb{Z}),\, \zeta \mapsto c_1(\zeta)$ of $H$-spaces,
where $c_1$ is the first Chern class).

Thus, after taking the limit as $k\rightarrow \infty$ in (\ref{obstruc})
we get the fibration
\begin{equation}
\label{usfib}
\begin{array}{ccccc}
\mathbb{C}P^{\infty} & \stackrel{i}{\hookrightarrow} &
\B{\cal U} & &  \mathbb{C}P^{\infty}
\hookrightarrow \prod_{q\geq 1}\K(\mathbb{Q},\, 2q)
\qquad \qquad \\
&& \downarrow{\scriptstyle p} && \downarrow \\
&& \B{\cal PU}, & \hbox{i.e.} &
\qquad \K(\mathbb{Q}/\mathbb{Z},\, 2)\times
\prod_{q\geq 2}\K(\mathbb{Q},\, 2q). \\
\end{array}
\end{equation}

The $H$-space $\B{\cal PU}$ represents
the homotopy functor $X\mapsto \widetilde{\AB}(X)$
on the category of finite $CW$-complexes.
Indeed, since we take the direct limit $\B{\cal PU}=
\varinjlim_k\BPU(k)$ over the maps
of classifying spaces $\BPU(k)\stackrel{i_{k,\, m}}
{\longrightarrow}\BPU(km)$ such that
$i^*_{k,\, m}(A^{univ}_{km})=
A_k^{univ}\otimes \widetilde{M}_m$ (where by
$A_r^{univ}$ we denote the
universal $M_r(\mathbb{C})$-bundle
over $\BPU(r)$), the stable
equivalence relation (\ref{classicequiv}) appears.

Now it is easy to see that $\Br(X)=\coker \{p_*\colon
[X,\, \B{\cal U}]\rightarrow
[X,\, \B{\cal PU}]\} ,$ i.e. $\Br(X)=\coker \{
[X,\, \K(\mathbb{Q},\, 2)]\rightarrow
[X,\,  \K(\mathbb{Q}/\mathbb{Z},\, 2)]\}=\coker \{
H^2(X;\, \mathbb{Q})\rightarrow H^2(X;\,
\mathbb{Q}/\mathbb{Z})\}.$
Using the exact sequence of cohomology groups
induced by the sequence of coefficients
$$
0\rightarrow \mathbb{Z}\rightarrow \mathbb{Q}\rightarrow
\mathbb{Q}/\mathbb{Z}\rightarrow 0,
$$
we see that
$\coker \{
H^2(X;\, \mathbb{Q})\rightarrow H^2(X;\, \mathbb{Q}/\mathbb{Z})\}=
\im \{H^2(X;\, \mathbb{Q}/\mathbb{Z})
\stackrel{\delta}{\rightarrow}
H^3(X;\, \mathbb{Z})\}$ (here $\delta$ is the coboundary
homomorphism), i.e. $\Br(X)=H^3_{tors}(X;\, \mathbb{Z})$
\cite{Grothendieck}.

The explicit form of the isomorphism $\Br(X)\cong
H^3_{tors}(X;\, \mathbb{Z})$ can be described as follows.
Recall that the structure group of a bundle $A_k$ with
fiber $M_k(\mathbb{C})$ is
$\Aut(M_k(\mathbb{C}))=\PGL_k(\mathbb{C})$ and $\PGL_k(\mathbb{C})$
contains $\PU(k)$ as a strong deformation retract.
The obstruction theory asserts (\cite{Fuchs}, Ch.5) that the first
(and unique!) obstruction
for the lifting of $\PU(k)$-bundle to
a $\U(k)$-bundle belongs to the group
$H^3(X;\, \pi_2(\mathbb{C}P^\infty))$
(see fibration (\ref{obstruc})). Therefore the assignment
$$
A_k\mapsto \; \{ \hbox{the first obstruction for the lifting}\}
$$
gives us the required description of the isomorphism
$\Br(X)\cong H^3_{tors}(X;\, \mathbb{Z}).$

Thus, any element of the group $H^3_{tors}(X;\, \mathbb{Z})$
can be realized as an obstruction for the lifting
of some $\PU$-bundle to a $\U$-bundle
(or equivalent a $\PGL$-bundle to a $\GL$-bundle).

\begin{remark}
Let us remark that it is not necessarily
that for $\alpha \in H^3(X;\, \mathbb{Z})$
such that $k\alpha=0$ there exists a bundle $A_k$ with fiber $M_k(\mathbb{C})$
whose invariant is equal to $\alpha$ (see \cite{Atiyah}, p.11).
\end{remark}

\begin{remark}
\label{upto1}
Now suppose $A_k$ has a lifting $\xi_k$ (i.e. $\End(\xi_k)=A_k$).
Then the obstruction theory says that $\xi_k$
is determined by $A_k$
up to taking the tensor product with
a line bundle $\zeta$ over $X$.
Indeed, it was already mentioned that $H^2(X;\, \pi_2(\mathbb{C}P^\infty))=
H^2(X;\, \mathbb{Z})$ is isomorphic to the group of
line bundles with respect to the tensor product;
from the other hand, for any line bundle $\zeta$
we have $\End(\xi_k \otimes \zeta)=\End(\xi_k)$.
\end{remark}
\begin{remark}
\label{kprimcomp}
Note that for any fixed integer $k>1$
one can develop the corresponding
theory of bundles with fibers of the form
$M_{k^n}(\mathbb{C})$ for arbitrary $n\in \mathbb{N}.$
In this way one can define the $k$-primary component $\Br_k(X)$
of the Brauer group as
$\coker \{ p_{k*}\colon [X,\, \BU(k^\infty)]\rightarrow [X,\, \BPU(k^\infty)]\}$.
\end{remark}

\subsection{An obstruction for an embedding of a
locally trivial algebra bundle into a trivial one}

In this subsection using the
$\Hom(M_k(\mathbb{C}),\, M_{kl}(\mathbb{C}))$-bundle
${\bf H}_{k,\, l}(A_k^{univ})\rightarrow \BPU(k)$
(which has been studied in Subsection 1.4)
we define topological obstructions for an embedding of a locally trivial
matrix algebra bundle into a trivial one. Note that we consider only
embeddings that are fiberwise homomorphisms of central algebras.

More precisely, let $A_k$ be an $M_k(\mathbb{C})$-bundle
over a finite $CW$-complex $X$.
In this subsection we construct a cohomological obstruction
for the existence of an embedding
$\mu \colon A_k\rightarrow X\times M_{kl}(\mathbb{C})$ such that
$\mu |_x (A_k)_x\subset M_{kl}(\mathbb{C})$
is a central subalgebra for any $x\in X$.
This obstruction equals $0$ iff $A_k$
is the core of some FAB $(A_k,\, \mu,\, \widetilde{M}_{kl})$ over $X$.

\begin{remark}
Note that nontrivial obstructions can exist (after taking the limit as in
Remark \ref{stabil} below) only if $(k,\, l)=1.$ In other words, all
the obstructions vanish in the stable case if we reject
this condition.
\end{remark}

Consider the fibration
\begin{equation}
\label{liftanz}
\begin{array}{c}
\diagram
\Fr_{k^n,\, l^n} \rto & {\bf H}_{k^n,\, l^n}(A_{k^n}^{univ}) \dto^{p_{k^n}}
& \!\!\!\!\!\!\!\!\!\!\!\! \simeq \Gr_{k^n,\, l^n} \\
& \BPU(k^n). & \\
\enddiagram
\end{array}
\end{equation}
Let $X$ be a finite $CW$-complex, $\dim(X)<2l^n$.
Then for a given map $f\colon X\rightarrow \BPU(k^n)$ we have the
first obstruction $\alpha(f)\in H^{2i}(X;\, \pi_{2i-1}(\Fr_{k^n,\, l^n}))$
for the lifting in (\ref{liftanz}) (\cite{Fuchs}, Ch.5),
where $\pi_{2i-1}(\Fr_{k^n,\, l^n})=\mathbb{Z}/k^n\mathbb{Z}$
and $\pi_{2i}(\Fr_{k^n,\, l^n})=0$
because $i<l^n.$ But we have shown in Subsection 1.4 that a lifting of
$f\colon X\rightarrow \BPU(k^n)$ in (\ref{liftanz})
is the same thing as an embedding $f^*(A_{k^n}^{univ})\hookrightarrow X\times M_{(kl)^n}(\mathbb{C})$.
Therefore the class $\alpha(f)$ is the
first obstruction for the existence of an embedding of
$f^*(A_{k^n}^{univ})$ into the
trivial bundle $X\times M_{(kl)^n}(\mathbb{C})$ (recall that $(k,\, l)=1$).

Note that for a mapping $f\colon X\rightarrow \BU(k^n)$
we can also define the obstruction for an embedding of the bundle
$f^*(\End(\xi_{k^n}^{univ}))=\End(f^*(\xi_{k^n}^{univ}))$
into the trivial
one $X\times M_{(kl)^n}(\mathbb{C})$. For this purpose instead of
bundle (\ref{liftanz}) we should consider the bundle
\begin{equation}
\label{liftanzu}
\begin{array}{c}
\diagram
\Fr_{k^n,\, l^n} \rto & {\bf H}_{k^n,\, l^n}
(\End({\xi}_{k^n}^{univ})) \dto^{\widehat{p}_{k^n}}
& \!\!\!\!\!\!\!\!\!\!\!\! \simeq \widehat{\Gr}_{k^n,\, l^n} \\
& \BU(k^n) & \\
\enddiagram
\end{array}
\end{equation}
(see Remark \ref{unitver}). Clearly,
the obstructions will be the same as in the ``projective'' case
(in connection with this note that the map
$\BU(k^n)\rightarrow \BPU(k^n)$ induces an
isomorphism of cohomology groups
$H^2(\BPU(k^n);\, \mathbb{Z}/k^n\mathbb{Z})
\cong H^2(\BU(k^n);\, \mathbb{Z}/k^n\mathbb{Z})$).

Now take $k,\, l,\, m,\, n$ such that $(km,\, ln)=1.$
Notice that there is the natural map
$$
\phi_{k,\,l;\, m,\,n}\colon \Hom(M_k(\mathbb{C}),\, M_{kl}(\mathbb{C}))\times
\Hom(M_m(\mathbb{C}),\, M_{mn}(\mathbb{C}))\rightarrow
\Hom(M_{km}(\mathbb{C}),\, M_{klmn}(\mathbb{C}))
$$
induced
by the tensor product of matrix algebras.

\begin{remark}
\label{Hspacestr}
If we identify $\Hom(M_k(\mathbb{C}),\, M_{kl}(\mathbb{C}))$
with $\Fr_{k,\, l}$ then the map $\phi_{k,\,l;\, m,\,n}$ can be described as
follows. One can easily verify that for a $k$-frame
$\{ \alpha_{i,\, j}\mid 1\leq i,\, j\leq k\}$ in $M_{kl}(\mathbb{C})$
and an m-frame $\{ \beta_{r,\, s}\mid 1\leq r,\, s\leq m\}$
in $M_{mn}(\mathbb{C})$ the collection
$\{ \alpha_{i,\, j}\otimes \beta_{r,\, s}\mid 1\leq i,\, j\leq k,\; 1
\leq r,\, s\leq m\}$ is a $km$-frame in $M_{klmn}(\mathbb{C})=
M_{kl}(\mathbb{C})\otimes M_{mn}(\mathbb{C})$.
Thus we have the natural map $\Fr_{k,\, l}\times \Fr_{m,\, n}
\rightarrow \Fr_{km,\, ln}$ which coincides with $\phi_{k,\,l;\, m,\,n}$
under the mentioned identification.
\end{remark}

Clearly, we also have the corresponding map of bundles
$$
\varphi_{k,\,l;\, m,\,n}\colon {\bf H}_{k,\, l}(A_k^{univ})\times
{\bf H}_{m,\, n}(A_m^{univ})\rightarrow
{\bf H}_{km,\, ln}(A_{km}^{univ})
$$
such that the diagram
$$
\diagram
{\bf H}_{k,\, l}(A_k^{univ})\times {\bf H}_{m,\, n}(A_m^{univ})
\rto^{\qquad \; \; \: \varphi_{k,\,l;\, m,\,n}} \dto_{p\times p} &
{\bf H}_{km,\, ln}(A_{km}^{univ}) \dto^p \\
\BPU(k)\times \BPU(m)\rto^{\qquad \otimes} & \BPU(km) \\
\enddiagram
$$
commutes.
Obviously, the map $\varphi_{k,\,l;\, m,\,n}$ can be identified with
the natural map $\Gr_{k,\, l}\times \Gr_{m,\, n}\rightarrow
\Gr_{km,\, ln}$ of matrix Grassmannians induced by the tensor product of
matrix algebras. These maps define the above described (see Subsection 2.2) $H$-space structure
on $\Gr:=\varinjlim_{(k,\, l)=1}\Gr_{k,\, l}\cong \BSU_{\otimes}.$

The same is true in the ``unitary'' case.

\begin{lemma}
\label{innd}
The homotopy type of $\Fr_{k,\, l^{\infty}}$ does not depend on
the choice of $l,\, (k,\, l)=1$.
\end{lemma}
{\noindent {\it Proof}.}\quad
Take $m$ such that $(k,\, m)=1.$ Define the map
$\alpha_{k,\, l,\, m}\colon \Fr_{k,\, l}\rightarrow \Fr_{k,\, lm}$
by the commutative diagram:
\begin{equation}
\label{innddd}
\diagram
\PU(lm)\rto^{E_{k}\otimes \ldots} &
\quad \PU(klm)\rto &
\Fr_{k,\, lm} \\
\PU(l)
\rto^{E_{k}\otimes \ldots} \uto^{\ldots \otimes E_m} &
\; \PU(kl)\rto
\uto_{\ldots \otimes E_m} &
\, \Fr_{k,\, l} \uto_{\alpha_{k,\, l,\, m}} \\
\enddiagram
\end{equation}
(the rows are fibrations).
One can easily verify that
the map $\alpha_{k,\, l,\, m}\colon \Fr_{k,\, l}\rightarrow \Fr_{k,\, lm}$
is a homotopy equivalence in dimensions $<2l$.
Now using the diagram
\begin{equation}
\label{yetanotherdiagr}
\diagram
& \Fr_{k,\, lm} \\
\Fr_{k,\, l} \urto^{\alpha_{k,\, l,\, m}} & &
\Fr_{k,\, m} \ulto_{\alpha_{k,\, m,\, l}} \\
\enddiagram
\end{equation}
we get the desired assertion$.\quad \square$

\smallskip

The proved lemma allows us to omit $l$ in the following notation.

Suppose $(m,n)=1,\: k|m,\, l|n$ ($\Rightarrow (k,l)=1$).
Consider the following morphism of fiber bundles:
\begin{equation}
\nonumber
\begin{array}{c}
\diagram
&\PU(n)\rto^{E_m\otimes\ldots} &
\PU(mn)\dto\\
\PU(l)\urto^{\ldots \otimes E_{\frac{n}{l}}}
\rto^{E_{k}\otimes\ldots} &
\PU(kl) \dto \urto^{E_{\frac{m}{k}}\otimes \ldots
\otimes E_{\frac{n}{l}}} &
\Fr_{m,\, n} \\
& \Fr_{k,\, l} \urto^\alpha
\enddiagram
\end{array}
\end{equation}
Using such maps $\alpha$ we can form the
direct limit $\Fr_{k^{\infty},\, l^{\infty}}$
whose homotopy type does not depend on the choice of $l,$
so we shall denote it just by $\Fr_{k^\infty}$.

\begin{lemma}
\label{homgrfr}
$\pi_r(\Fr_{k^\infty})=\varinjlim_{n}
\mathbb{Z}/k^n\mathbb{Z}$ for odd $r$ and $0$ otherwise.
Moreover, the natural embedding
$\Fr_{k^n}\hookrightarrow \Fr_{k^\infty}$
induces the monomorphisms $\pi_{2r-1}(\Fr_{k^n})\cong \mathbb{Z}/k^n\mathbb{Z}
\hookrightarrow \pi_{2r-1}(\Fr_{k^\infty}),\: 1\leq r\leq \infty$.
\end{lemma}
{\noindent {\it Proof}}\; follows from simple calculations with homotopy sequences
of obvious fibrations.$\quad \square$

\smallskip

So we can define the direct limit (as $n\rightarrow \infty$) of (\ref{liftanz})
\begin{equation}
\label{liftanz2}
\begin{array}{c}
\diagram
\Fr_{k^\infty} \rto & {\bf H}_{k^\infty}(A_{k^\infty}^{univ}) \dto^{p_k}
& \!\!\!\!\!\!\!\!\!\!\!\! \simeq \Gr \\
& \BPU(k^\infty) & \\
\enddiagram
\end{array}
\end{equation}
and of (\ref{liftanzu})
\begin{equation}
\label{liftanz3}
\begin{array}{c}
\diagram
\Fr_{k^\infty} \rto & {\bf H}_{k^\infty}(\End(\xi_{k^\infty}^{univ})) \dto^{\widehat{p}_k}
& \!\!\!\!\!\!\!\!\!\!\!\! \simeq \widehat{\Gr} \\
& \BU(k^\infty), & \\
\enddiagram
\end{array}
\end{equation}
respectively. All the mappings are the homomorphisms of $H$-spaces. In particular, we see that
$p_k,\, \widehat{p}_k$ are fiber substitutes for maps $t_k,\, \widehat{t}_k$ (see Section 2), respectively.

\begin{remark}
\label{stabil}
Suppose we are given a map $f\colon X\rightarrow \BPU(k^n)$ such that
$\alpha(f)\neq 0,\; \alpha(f)\in H^{2i}(X;\,\pi_{2i-1}(\Fr_{k^n,\, l^n}))$
(see fibration (\ref{liftanz})), $\dim X<2l^n$. Then the corresponding limit
obstruction in (\ref{liftanz2}) is also nontrivial. One can prove this using the following fact:
for $(km,\, ln)=1$ the natural map $\Fr_{k,\, l}\rightarrow \Fr_{km,\, ln}$
induces monomorphisms of the homotopy groups:
$\mathbb{Z}/k\mathbb{Z}\cong \pi_{2i-1}(\Fr_{k,\, l})
\hookrightarrow \pi_{2i-1}(\Fr_{km,\, ln})\cong \mathbb{Z}/km\mathbb{Z},\;
i\leq l.$
\end{remark}

\begin{remark}
\label{someeqrel}
Let $\F_k(X)$ be the set of equivalence classes of
$M_{k^n}(\mathbb{C})$-bundles $A_{k^n}$ (for arbitrary $n\in \mathbb{N}$)
with respect to the following equivalence relation:
$$
A_{k^m}\sim A_{k^n} \Leftrightarrow \hbox{ there are FAB's cores }
B_{k^r},\, B_{k^s}
$$
$$
\hbox{ such that } A_{k^m}\otimes B_{k^r}
\cong A_{k^n}\otimes B_{k^s} \; (\Rightarrow m+r=n+s).
$$
Obviously, $\F_k(X)$ is an Abelian group with respect to the
operation induced by the tensor product of bundles.
Moreover, $\F_k(X)=\coker \{ p_{k*}\colon [X,\, \Gr]
\rightarrow[X,\, \BPU(k^\infty)]\}$ (=$\im \{ h_{k*}\colon [X,\, \BPU(k^\infty)]
\rightarrow [X,\, \BN_{k^\infty}^\times]\}$, see the next subsection).
Note that the first cohomological obstruction for the lifting in (\ref{liftanz2})
is well-defined on such equivalence classes.

Although the group $\F_k(X)$ looks like the Brauer group $\Br_k(X)$, the actual
analog of the latter group will be defined in the next subsection.
\end{remark}

\subsection{A generalized Brauer group}

Consider exact sequence (\ref{liftanz3}) of $H$-spaces. Our goal is to extend it to the right.
So at the next step we have to define a ``classifying space'' for $\Fr_{k^\infty}.$
One possible approach uses the fact that $\Fr_{k^\infty}$ is an infinite loop space
(as a fiber of the localization map
$\widehat{t}_{k}\colon
\BU_{\otimes}=\widehat{\Gr}\rightarrow \BU[\frac{1}{k}]_{\otimes}
=\BU(k^{\infty})$ for $\BU_{\otimes}$; $\BU_{\otimes}$ is an infinite loop
space due to G.B. Segal \cite{Segal}).
In this way one obtains the fibration
\begin{equation}
\label{extended2}
\widehat{\Gr}\stackrel{\widehat{t}_k}{\longrightarrow}
\BU(k^\infty)\stackrel{\widehat{h}_{k}}
{\longrightarrow}\BFr_{k^\infty},
\end{equation}
where $\BFr_{k^\infty}$ is the base of the universal principal
$\Fr_{k^\infty}$-fibration, $\Omega(\BFr_{k^\infty})=\Fr_{k^\infty}.$

But we present a more direct approach by replacing the loop spaces
by groups of the same homotopy type. More precisely,
we introduce topological groups $\GL_{k^\infty}({\cal K}(H)),\, \N^{\times}_{k^{\infty}},\,
\GL_1(\Delta_{k^\infty})^0=\GL_1({\cal C}(H))^0$ together with homomorphisms
$\GL_{k^\infty}({\cal K}(H))\rightarrow \N^{\times}_{k^{\infty}},\,
\N^{\times}_{k^{\infty}}\rightarrow \GL_1(\Delta_{k^\infty})^0$
forming the exact sequence of groups
$$
\GL_{k^\infty}({\cal K}(H))\rightarrow \N^{\times}_{k^{\infty}}\rightarrow \GL_1(\Delta_{k^\infty})^0
$$
which is homotopy equivalent to the fibration
$$
\U(k^\infty)\rightarrow \Fr_{k^\infty}\rightarrow \widehat{\Gr}
$$
(see diagram (\ref{unitver1})). Then we identify
$\BFr_{k^\infty}$ with the classifying space
$\BN^{\times}_{k^{\infty}}$ and the map $\widehat{h}_k\colon
\BU(k^\infty)\rightarrow \BFr_{k^\infty}$ with the map of
classifying spaces $\BGL_{k^\infty}({\cal K}(H))\rightarrow
\BN^{\times}_{k^{\infty}},$ induced by the above group
homomorphism. The main advantage of this approach is that any element of
the defined below generalized Brauer group can be represented by a locally
trivial bundle with the structure group $\N_k^\times$ (for some $k\in \mathbb{N}$).

So let $H$ be a separable Hilbert space, ${\cal B}(H)$ and ${\cal K}(H)$
the algebra of bounded operators in $\End(H)$ and the ideal of
compact operators in ${\cal B}(H)$, respectively.
By ${\cal C}(H)$ denote the Calkin algebra ${\cal B}(H)/{\cal K}(H).$
Put
$$
\Delta_k:=\left\{ \begin{pmatrix}
\lambda & 0 & \ldots & 0 \\
0 & \lambda & \dots & 0 \\
\vdots & \vdots & \ddots & \vdots \\
0 & 0 & \ldots & \lambda \\
\end{pmatrix} \mid \lambda \in {\cal C}(H) \right\}\subset M_k({\cal C}(H)).
$$
Clearly, $\Delta_k$ is a subalgebra in $M_k({\cal C}(H))$ isomorphic to ${\cal C}(H).$
Let $\pi_k \colon M_k({\cal B}(H))\rightarrow M_k({\cal C}(H))$ be the natural epimorphism.
Let $\N_k$ be the subalgebra $\pi_k^{-1}(\Delta_k)\subset M_k({\cal B}(H))$,
$\N_k^\times:=\GL_1(\N_k)=\N_k \cap \GL_k({\cal B}(H))$ its multiplicative group.
It is a closed subgroup (in the norm topology) in $\GL_k({\cal B}(H))$. The groups $\N_{k^n}^\times$
play a crucial role in the further constructions. In particular, we shall show that in contrast to
the group $\GL_1({\cal B}(H))$, they are not contractible, if $n>0$
(note that $\N_{k^0}^\times \cong \GL_1({\cal B}(H))$).

\begin{remark}
One can replace $\N_k^\times$ by the homotopy equivalent (polar decomposition!)
unitary group $\N_k^\times \cap \U_k({\cal B}(H)),$ where $\U_k({\cal B}(H))\subset M_k({\cal B}(H))$
is the group of unitary operators contained in $M_k({\cal B}(H))$.
\end{remark}

\begin{remark}
Let us remark that there is an obvious way to transfer our definition of
$\N^\times_k$ to the algebraic $K$-theory of a ring $R$
(using the analogy between ${\cal C}(H)$ and $\Sigma R$).
\end{remark}

The following results about Calkin algebra are needed for the sequel.
Let $\pi \colon {\cal B}(H)\rightarrow {\cal C}(H)$ be the canonical epimorphism.
By $\GL_1({\cal C}(H))^0$ denote the connected component of the unit in $\GL_1({\cal C}(H)).$
In other words, it is the image of the space $\Fred(H)^0$ of zero index Fredholm operators
under the map $\pi |_{\Fred(H)}\colon \Fred(H)\rightarrow \GL_1({\cal C}(H))$
(which is a fibration with fiber an affine space over ${\cal K}(H)$).
Let $\bar{\pi}\colon \GL_1({\cal B}(H))\rightarrow \GL_1({\cal C}(H))^0$ be the birestriction of
$\pi$ to the multiplicative
subgroups. Clearly, $\ker(\bar{\pi})=\GL_1({\cal K}(H)),$ where $\GL_1({\cal K}(H))$ is the group
of invertible operators of the form $1+K,\; K\in {\cal K}(H).$ It is known \cite{Palais}
that this
group is homotopy equivalent to the infinite unitary group $\U:=\varinjlim_n\U(n)$
(with respect to the standard inclusions).
Recall also that the group $\GL_1({\cal C}(H))$ ($\GL_1({\cal C}(H))^0$) is homotopy equivalent to
$\mathbb{Z}\times \BU$ ($\BU$ respectively).
Finally note that $\GL_k({\cal K}(H))
:=\{ \hbox{the group of invertible
$k\times k$-matrices of the form } 1+K \mid K \in M_k({\cal K}(H))\}$
is isomorphic to $\GL_1({\cal K}(H)).$

Consider the group epimorphisms $\bar{\pi}_k\colon \GL_k({\cal B}(H))\rightarrow
\GL_k({\cal C}(H))^0$ and $\widetilde{\pi}_k\colon
\N_k^\times \rightarrow \GL_1(\Delta_k)^0$ induced by $\pi_k$.
Clearly their kernels coincide with the subgroup $\GL_k({\cal
K}(H))$.

\begin{lemma}
\label{homgrcal}
The topological group $\N_k^\times$ has the following homotopy groups: $\pi_{2r}(\N_k^\times)=0,\;
\pi_{2r-1}(\N_k^\times)=\mathbb{Z}/k\mathbb{Z},\; r\geq 1$.
\end{lemma}
{\noindent \it Proof}\; easily follows from the commutative diagram of exact sequences
\begin{equation}
\label{diagrhgr}
\diagram
1 \rto & \GL_k({\cal K}(H)) \rto & \GL_k({\cal B}(H)) \rto^{\bar{\pi}_k} & \GL_k({\cal C}(H))^0 \rto & 1 \\
1 \rto & \GL_k({\cal K}(H)) \uto^= \rto & \N_k^\times \rto^{\widetilde{\pi}_k} \uto^\cup & \GL_1(\Delta_k)^0
\uto^\cup \rto & 1.
\enddiagram
\end{equation}
Indeed, recall that the group $\GL_k({\cal B}(H))\cong
\GL_1({\cal B}(H))$ is contractible (the polar decomposition + Kuiper's theorem).
In particular, upper row is just the path fibration.
Now the required assertion follows from the nontrivial piece
of the corresponding homotopy sequences:
$$
\diagram
0\rto & 0 \rto & \pi_{2r}(\GL_k({\cal C}(H))^0)\rto^{\cong \quad} &
\pi_{2r-1}(\GL_k({\cal K}(H))) \rto & 0 \rto & 0 \\
0 \rto & \pi_{2r}(\N_k^\times) \uto \rto & \pi_{2r}(\GL_1(\Delta_k)^0) \rto^{\cdot k \quad}
\uto^{\cdot k} &
\pi_{2r-1}(\GL_k({\cal K}(H))) \rto \uto^= & \pi_{2r-1}(\N_k^\times) \uto \rto & 0, \\
\enddiagram
$$
where $\cdot k$ means the homomorphism which takes the group generator $1$ to
$k\cdot 1.\quad \square$

\begin{remark}
\label{addrem} Since the inclusion
$\GL_1(\Delta_k)^0\hookrightarrow \GL_k({\cal C}(H))^0$ in diagram
(\ref{diagrhgr}) is a classifying map for the lower row
(considered as a $\GL_k({\cal K}(H))$-fibration), we have the fibration
$$
\N_k^\times \stackrel{\widetilde{\pi}_k}{\longrightarrow} \GL_1(\Delta_k)^0\longrightarrow \GL_k({\cal C}(H))^0.
$$
Because of $\Omega(\GL_1(\Delta_k)^0)\simeq \GL_1({\cal K}(H)),\:
\Omega(\GL_k({\cal C}(H))^0)\simeq \GL_k({\cal K}(H)),$
we also obtain the fibration
\begin{equation}
\label{extended30}
\diagram
\GL_1({\cal K}(H)) \rto & \GL_k({\cal K}(H)) \rto & \N_k^\times,
\enddiagram
\end{equation}
where the inclusion $\GL_1({\cal K}(H))\rightarrow \GL_k({\cal K}(H))$ is the following:
\begin{equation}
\label{avermat}
\alpha \mapsto \begin{pmatrix}
\alpha & 0 & \ldots & 0 \\
0 & \alpha & \dots & 0 \\
\vdots & \vdots & \ddots & \vdots \\
0 & 0 & \ldots & \alpha \\
\end{pmatrix}
\end{equation}
($k\times k$-matrix); in particular, $\N_k^\times$ is homotopy equivalent to
the homogeneous space $\GL_k({\cal K}(H))/\GL_1({\cal K}(H)).$
\end{remark}

\begin{remark}
\label{addremm}
Since $\BGL_1({\cal K}(H))\simeq \GL_1(\Delta_{k})^0$ we see that there is the map
$\N_{k}^\times \rightarrow \BGL_1({\cal K}(H))$ which is equivalent to the epimorphism
$\widetilde{\pi}_{k}\colon \N_{k}^\times \rightarrow \GL_1(\Delta_{k})^0,\;
\ker(\widetilde{\pi}_{k})=\GL_{k}({\cal K}(H))$; analogously,
there is the map
$\BGL_1({\cal K}(H))\rightarrow \BGL_{k}({\cal K}(H))$ which is equivalent to the natural
inclusion $\GL_1(\Delta_{k})^0\rightarrow \GL_{k}({\cal C}(H))^0,$
where ${\cal C}(H)$ is the Calkin algebra as above.
\end{remark}

Consider the homomorphisms $\U(kn)\stackrel{\kappa_n}{\hookrightarrow}\U(k(n+1)),$
$$
\begin{pmatrix}
\alpha_{11} & \alpha_{12} & \ldots & \alpha_{1k} \\
\alpha_{21} & \alpha_{22} & \dots & \alpha_{2k} \\
\vdots & \vdots & \ddots & \vdots \\
\alpha_{k1} & \alpha_{k2} & \ldots & \alpha_{kk} \\
\end{pmatrix}
\mapsto
\begin{pmatrix}
\alpha_{11} & 0 & \alpha_{12} & 0 & \ldots & \alpha_{1k} & 0 \\
0 & 1 & 0 & 0 & \ldots & 0 & 0 \\
\alpha_{21} & 0 & \alpha_{22} & 0 & \dots & \alpha_{2k} & 0 \\
0 & 0 & 0 & 1 & \ldots & 0 & 0 \\
\vdots & \vdots & \vdots & \vdots & \ddots & \vdots & \vdots \\
\alpha_{k1} & 0 & \alpha_{k2} & 0 & \ldots & \alpha_{kk} & 0 \\
0 & 0 & 0 & 0 & \ldots & 0 & 1 \\
\end{pmatrix},
$$
where $\alpha_{ij}\in M_{n}(\mathbb{C}).$ Let $i_n\colon \U(n)\hookrightarrow \U(n+1)$
be the standard inclusion
$$
\alpha \mapsto \begin{pmatrix}
\alpha & 0 \\
0 & 1 \\
\end{pmatrix},\; \alpha \in \U(n).
$$
Consider also the homomorphisms
$\U(n)\stackrel{E_k\otimes \ldots}{\longrightarrow}\U(kn)$. One can easily check
that $E_k\otimes (i_n\alpha)=\kappa_n(E_k\otimes \alpha)\: \forall \alpha \in \U(n).$ Thus, we obtain
the well-defined map of the direct limits:
\begin{equation}
\label{tensmapp}
\tau_k\colon \U =\varinjlim_{i_n}\U(n)\rightarrow \varinjlim_{\kappa_n}\U(kn)\cong \U
\end{equation}
induced by the tensor product with $E_k$ (cf. diagram (\ref{avermat})).

\begin{lemma}
\label{klrespprime}
The space $\Fr_{k}$ is the fiber of the map $\B \tau_k\colon \BU \rightarrow \BU$
induced by the group homomorphism (\ref{tensmapp}).
\end{lemma}
{\noindent \it Proof.}\; Obviously, the map
$\BU(l^n)\stackrel{[k]\otimes \ldots}{\longrightarrow}\BU(kl^n)$ induced by the tensor product
with a trivial bundle of rank $k$ is equivalent to the map of these classifying spaces induced
by the group homomorphism $\U(l^n)\stackrel{E_k\otimes \ldots}{\longrightarrow}\U(kl^n).$
Hence $\Fr_{k,\, l^n}$ is the fiber of
$\BU(l^n)\stackrel{[k]\otimes \ldots}{\longrightarrow}\BU(kl^n).$
It easily follows from the commutative
diagram (cf. diagram (\ref{innddd}))
$$
\diagram
\Fr_{k,\, l^{n+1}} \rto & \BU(l^{n+1}) \rto^{[k]\otimes \ldots} & \BU(kl^{n+1}) \\
\Fr_{k,\, l^n} \uto \rto & \BU(l^n) \uto^{\ldots \otimes [l]} \rto^{[k]\otimes \ldots} &
\BU(kl^n) \uto^{\ldots \otimes [l]} \\
\enddiagram
$$
that the natural inclusions $\Fr_{k,\, l^n}\rightarrow \Fr_{k,\, l^{n+1}}$ are weak homotopy equivalences
($\Rightarrow$ homotopy equivalences because $\Fr_{k,\, l}$ is a $CW$-complex)
up to dimension $\sim 2l^n$ (recall that $(k,\, l)=1$).$\quad \square$

\smallskip

The following proposition makes the statement of Lemma \ref{homgrcal} more precise (cf. Lemma \ref{homgrfr}).

\begin{proposition}
\label{joinedprop} There is a homotopy equivalence $\Fr_{k}\simeq
\N^\times_k$ such that the diagram (whose rows are the above
fibrations (cf. Remarks \ref{addrem} and \ref{addremm}) and the
vertical arrows are homotopy equivalences)
$$
\diagram
\N_{k}^\times \rto & \BGL_1({\cal K}(H)) \rto & \BGL_{k}({\cal K}(H)) \\
\Fr_{k} \rto \uto^\simeq & \BU \rto^{\B \tau_k} \uto^\simeq & \BU \uto^\simeq \\
\enddiagram
$$
commutes up to homotopy.
\end{proposition}
{\noindent {\it Proof}.}\quad
We have the group homomorphism $\U \rightarrow \GL_1({\cal K}(H))$ which
is a homotopy equivalence. Clearly, this equivalence identifies the group homomorphisms
$\tau_k \colon \U \rightarrow \U$ and
$\GL_1({\cal K}(H))\rightarrow \GL_{k}({\cal K}(H))$ (see (\ref{avermat})).
There is also the homotopy equivalence
$\BU\rightarrow \BGL_1({\cal K}(H))$ of classifying spaces and
the corresponding maps of classifying spaces $\B \tau_k \colon \BU \rightarrow \BU$ and
$\BGL_1({\cal K}(H))\rightarrow \BGL_{k}({\cal K}(H))$ can also be identified.
Therefore their homotopy fibers $\Fr_{k}$ and $\N_{k}^\times$ are homotopy equivalent too, and
the diagram is commutative.$\quad \square$

\begin{remark}
Let us remark that the fibration (see (\ref{extended30}))
\begin{equation}
\label{extended33}
\BGL_1({\cal K}(H))\rightarrow \BGL_k({\cal K}(H))\rightarrow \BN_{k}^\times
\end{equation}
has the following interpretation. Let
$\EN_{k}^\times \rightarrow \BN_{k}^\times$ be the universal principal
$\N_{k}^\times$-bundle. Using the action of the subgroup $\GL_{k}({\cal K}(H))\subset
\N_{k}^\times$
on its fibers $\cong \N_{k}^\times,$ we obtain (\ref{extended33}) (because $\BGL_1({\cal K}(H))\cong
\GL_1(\Delta_{k})^0\cong
\N_{k}^\times / \GL_{k}({\cal K}(H))$).
\end{remark}

Now consider the commutative diagram of group homomorphisms:
$$
\diagram
\GL_{l^{n+1}}({\cal K}(H)) \rto^{E_k\otimes \ldots} & \GL_{kl^{n+1}}({\cal K}(H)) \rto & \GL_{l^{n+1}}(\N_k) \\
\GL_{l^n}({\cal K}(H)) \rto^{E_k\otimes \ldots} \uto^{\ldots \otimes E_l} & \GL_{kl^n}({\cal K}(H))
\rto \uto^{\ldots \otimes E_l} & \GL_{l^n}(\N_k). \uto^{\ldots \otimes E_l} \\
\enddiagram
$$
Since $(k,\, l)=1,$ we see that the arrow $\GL_{l^n}(\N_k)\rightarrow \GL_{l^{n+1}}(\N_k)$
is a homotopy equivalence (cf. the proof of Lemma \ref{klrespprime}).
Put $\GL_{l^\infty}({\cal K}(H)):=\varinjlim_n\GL_{l^n}({\cal K}(H)),\;
\GL_{kl^\infty}({\cal K}(H)):=\varinjlim_n\GL_{kl^n}({\cal K}(H))$, where the direct limits are taken over the group
homomorphisms $\GL_{l^n}({\cal K}(H))\stackrel{\ldots \otimes E_l}{\longrightarrow}\GL_{l^{n+1}}({\cal K}(H))$ and
$\GL_{kl^n}({\cal K}(H))\stackrel{\ldots \otimes E_l}{\longrightarrow}\GL_{kl^{n+1}}({\cal K}(H))$
($\GL_{l^\infty}(\N_k)\simeq \GL_1(\N_k)=\N_k^\times$). Thus, we have the fibration
$$
\GL_{l^\infty}({\cal K}(H))\rightarrow \GL_{kl^\infty}({\cal K}(H))\rightarrow \N_k^\times
$$
and hence also the fibration
$$
\N_k^\times \rightarrow \BGL_{l^\infty}({\cal K}(H))\stackrel{[k]\otimes \ldots}
{\longrightarrow}\BGL_{kl^\infty}({\cal K}(H)).
$$

Taking into account the obtained results, the proof of the following lemma is trivial.
\begin{lemma}
\label{relathomfunc}
There is the commutative (up to homotopy) diagram
\begin{equation}
\label{relathomfunc1}
\diagram
\N_k^\times \rto & \BGL_{l^\infty}({\cal K}(H))\rto^{[k]\otimes \ldots} & \BGL_{kl^\infty}({\cal K}(H)) \\
\Fr_k \uto^{\simeq} \rto & \BU(l^\infty) \uto^{\simeq} \rto^{[k]\otimes \ldots} & \BU(kl^\infty) \uto_{\simeq} \\
\enddiagram
\end{equation}
whose rows are fibrations and vertical arrows are homotopy equivalences.
\end{lemma}

\begin{lemma}
There is the commutative (up to homotopy) diagram
$$
\diagram
\Fr_{k,\, l^n} \rto & \BPU(l^n) \rto^{[k]\otimes \ldots} & \BPU(kl^n) \\
\Fr_{k,\, l^n} \uto^= \rto & \Gr_{k,\, l^n} \uto \rto & \BPU(k) \uto_{\ldots \otimes [l^n]} \\
\enddiagram
$$
whose rows are the above fibrations.
\end{lemma}
{\noindent \it Proof.}\; Consider the mapping
$$
\diagram
{\rm Prin}(A_k^{univ}\otimes \widetilde{M}_{l^n}) \rto \dto & {\rm Prin}(A_{kl^n}^{univ}) \dto \\
\BPU(k) \rto^{\ldots \otimes [l^n]} & \BPU(kl^n),
\enddiagram
$$
of principal $\PU(kl^n)$-bundles (so it is $\PU(kl^n)$-equivariant).
We have the free action of the subgroup $E_k\otimes \PU(l^n)\subset \PU(kl^n)$ on their total
spaces. Note that ${\rm Prin}(A_k^{univ}\otimes \widetilde{M}_{l^n})$ is the fiber of the map
$\BPU(k)\stackrel{\ldots \otimes [l^n]}{\longrightarrow}\BPU(kl^n),$ hence
$$
{\rm Prin}(A_k^{univ}\otimes \widetilde{M}_{l^n})\cong
\PU(kl^n)/(\PU(k)\otimes E_{l^n}).
$$
Therefore its factor by $E_k\otimes \PU(l^n)$ is $\Gr_{k,\, l^n}$ and we obtain the required
commutative diagram.$\quad \square$

\smallskip

Replacing $\BPU$ by $\BU$ and $\Gr$ by $\widehat{\Gr}$ one obtains the analogous diagram
in the unitary case.

Put $\widehat{\Gr}_k:=\varinjlim_n\widehat{\Gr}_{k,\, l^n}.$ We claim that this homotopy type does not depend
on the choice of $l,\, (k,\, l)=1.$ Indeed, the maps $\alpha_{k,\, l,\, m}$ and $\alpha_{k,\, m,\, l}$ in diagram
(\ref{yetanotherdiagr}) are $\PU(k)$-equivariant.

\begin{corollary}
There is the commutative (up to homotopy) diagram
$$
\diagram
\Fr_k \rto & \BU(l^\infty) \rto^{[k]\otimes \ldots} & \BU(kl^\infty) \\
\Fr_k \uto^= \rto & \widehat{\Gr}_k \uto \rto & \BU(k) \uto \\
\enddiagram
$$
whose rows are the fibrations.
\end{corollary}

Combining the previous corollary with Lemma \ref{relathomfunc}, we obtain the following result.

\begin{theorem}
\label{thimportfib}
There is the fibration
\begin{equation}
\label{importfib}
\widehat{\Gr}_k \rightarrow \BU(k) \rightarrow \BN_k^\times,
\end{equation}
which actually is an extension of $\Fr_k\rightarrow \widehat{\Gr}_k\rightarrow \BU(k)$ to the right.
\end{theorem}
{\noindent \it Proof.}\; The upper fibration in diagram (\ref{relathomfunc1}) can be extended to the right
up to fibration
$$
\BGL_{l^\infty}({\cal K}(H))\rightarrow \BGL_{kl^\infty}({\cal K}(H))\rightarrow \BN_k^\times.
$$
The map $\BU(k) \rightarrow \BN_k^\times$ is just the composition $\BU(k)\rightarrow \BU(kl^\infty)
\simeq \BGL_{kl^\infty}({\cal K}(H))\rightarrow \BN_k^\times.$ It follows from the previous corollary
that the homotopic fiber of this map is just $\widehat{\Gr}_k. \quad \square$

\begin{remark}
Let us note an analogy between fibrations (\ref{obstruc}) and (\ref{importfib}) (cf. the definition of
generalized Brauer group below).
\end{remark}

There is a continuous action of the group $\N_k^\times$ on different spaces; for example one can consider the action
on the algebra $\N_k$ by inner automorphisms. So one can think of a bundle with the structure group
$\N_k^\times$ as a locally-trivial $\N_k$-bundle.

Let ${\bold F}_k$ be the functor on bundles, induced by the above map of classifying spaces
$\BU(k)\rightarrow \BN_k^\times$.
By analysis of the previous constructions, one can obtain the explicit form of the corresponding group
homomorphism $\U(k)\rightarrow \N_k^\times$. For example (for $k=2$) it is as follows:
\begin{equation}
\nonumber
\begin{pmatrix}
\alpha & \beta \\
\gamma & \delta \\
\end{pmatrix} \mapsto \begin{pmatrix}
\alpha & 0 & 0 & \ldots & \beta & 0 & 0 & \ldots \\
0 & 1 & 0 & \ldots & 0 & 0 & 0 & \ldots \\
0 & 0 & 1 & \ldots & 0 & 0 & 0 & \ldots \\
\vdots & \vdots & \vdots & \ddots & \vdots & \vdots & \vdots & \ddots \\
\gamma & 0 & 0 & \ldots & \delta & 0 & 0 & \ldots \\
0 & 0 & 0 & \ldots & 0 & 1 & 0 & \ldots \\
0 & 0 & 0 & \ldots & 0 & 0 & 1 & \ldots \\
\vdots & \vdots & \vdots & \ddots & \vdots & \vdots & \vdots & \ddots \\
\end{pmatrix}.
\end{equation}

\begin{corollary}
Suppose we are given an $\mathbb{C}^k$-bundle $\xi_k$ over a finite $CW$-complex $X$.
Then ${\bold F}_k(\xi_k)$ is a trivial
bundle iff $\End(\xi_k)$ is the core (see Subsection 2.3) of some FAB over $X$ (i.e. iff
there exists an embedding
$\End(\xi_k)\hookrightarrow X\times M_{km}(\mathbb{C})$ for some sufficiently
large $m,\, (k,\, m)=1$).
\end{corollary}
{\noindent {\it Proof.}\;} Since (\ref{importfib}) is a fibration, we have the corresponding exact sequence
of pointed sets $[X,\, \widehat{\Gr}_k]\rightarrow [X,\, \BU(k)]\rightarrow [X,\, \BN_k^\times].$ Now the required
assertion follows from the discussion after diagram (\ref{uncasedia}).$\quad \square$

\begin{remark}
Recall that we have already interpreted the lifting problem in the fibration $\Fr_k \rightarrow \widehat{\Gr}_k\simeq
{\bf H}_k(\End (\xi_k^{univ}))\rightarrow \BU(k)$ in terms of the existence of an embedding
into a trivial bundle, see the discussion around diagram (\ref{liftanzu}).
\end{remark}

\begin{remark}
Let us return to the fibration (\ref{importfib}). Since $\widehat{\Gr}_k=
\Gr_k{\mathop{\times}\limits_{\BPU(k)}}\BU(k)$ (see Subsection 2.6) it seems very reliable that the above map
$\BU(k)\rightarrow \BN_k^\times$ factors through $\BPU(k)$ i.e.
there exists a map $\BPU(k)\rightarrow \BN_k^\times$ (with fiber $\Gr_k$) such that
the diagram
$$
\diagram
\BU(k) \dto \drto & \\
\BPU(k) \rto & \BN_k^\times \\
\enddiagram
$$
commutes. This might show that the generalized Brauer group actually is a generalization of the classical
Brauer group (cf. diagram (\ref{actgen})).
\end{remark}

\begin{remark}
It is interesting to apply the obtained results to the equivalence relation, considered in Remark \ref{someeqrel}.
More precisely, we see that (in notation of the remark) $A_{k^m}\sim
A_{k^n} \Leftrightarrow {\bold F}_{k^m}(A_{k^m})$ and ${\bold F}_{k^n}(A_{k^n})$ are ``stable isomorphic''.
\end{remark}

It is obvious that the Kronecker product with the unit $k\times k$-matrix $E_k$ induces
the group homomorphisms $\N_{k^0}^{\times}\rightarrow \N_{k}^{\times}\rightarrow
\N_{k^{2}}^{\times}\rightarrow \ldots \;$.
For the topological group $\N_{k^n}^{\times}$ we have its classifying space $\BN_{k^n}^{\times}$.
Since $\B$ is a functor, we see that the above homomorphisms induce the corresponding
maps of classifying spaces $\BN_{k^0}^{\times}\rightarrow
\BN_{k}^{\times}\rightarrow \BN_{k^{2}}^{\times}\rightarrow \ldots \:$.

Consider the homomorphisms $\GL_1({\cal K}(H))\rightarrow \GL_k({\cal K}(H))\rightarrow \GL_{k^2}({\cal K}(H))
\rightarrow \ldots$ induced by the Kronecker product with the unit $k\times k$-matrix. Their direct
limit is the localization map $\GL_1({\cal K}(H))\rightarrow \GL_{k^\infty}({\cal K}(H))$.

\begin{corollary}
There is the diagram
$$
\diagram
\N_{k^\infty}^\times \rto & \BGL_1({\cal K}(H)) \rto & \BGL_{k^\infty}({\cal K}(H)) \\
\Fr_{k^\infty} \rto \uto^\simeq & \BU =\widehat{\Gr} \rto \uto^\simeq & \BU(k^\infty) \uto^\simeq \\
\enddiagram
$$
which is commutative up to homotopy.
\end{corollary}
{\noindent {\it Proof}\;} follows from Proposition \ref{joinedprop}.$\quad
\square$

\smallskip

Thus, instead of fibration (\ref{extended2})
we can consider the homotopy equivalent fibration
\begin{equation}
\label{extended3}
\BGL_1({\cal K}(H))\stackrel{\widehat{t}^\prime_k}{\longrightarrow} \BGL_{k^\infty}({\cal K}(H))
\stackrel{\widehat{h}^\prime_{k}}{\longrightarrow}\BN_{k^\infty}^\times
\end{equation}
(which is the direct limit of (\ref{extended33})).

Now we want to define a structure of homotopy commutative $H$-space on $\BN_{k^\infty}^\times$
which turn the above fibration into an exact sequence of $H$-spaces (i.e.
$$
[X,\, \BGL_1({\cal K}(H))]\stackrel{\widehat{t}^\prime_{k*}}{\longrightarrow}[X,\, \BGL_{k^\infty}({\cal K}(H))]
\stackrel{\widehat{h}^\prime_{k*}}{\longrightarrow}[X,\, \BN_{k^\infty}^\times]
$$
is an exact sequence of Abelian groups for any finite $CW$-complex $X$).
Unfortunately, we do not know any direct way to do this in terms of $\N_{k^n}^\times$-bundles.

\begin{remark}
For example, there are mappings
$\N_{k^m}^{\times}\times \N_{k^n}^{\times}\rightarrow \N_{k^{m+n}}^{\times}$
defined by the Kronecker product of matrices $M_{k^m}({\cal B}(H))\times M_{k^n}({\cal B}(H))
\rightarrow M_{k^{m+n}}({\cal B}(H))$ (cf. Remark \ref{Hspacestr}).
Indeed, one can easily verify (using the fact
that ${\cal K}(H)$  is a two-sided ideal in ${\cal B}(H)$) that
for any $A\in \N_{k^m}^{\times},\; B\in \N_{k^n}^{\times}$ their Kronecker product $A\otimes B$
satisfies the condition $\bar{\pi}_{k^{m+n}}(A\otimes B)\in \GL_1(\Delta_{k^{m+n}})^0\subset
\GL_{k^{m+n}}({\cal C}(H))^0,$
i.e. $A\otimes B \in \N_{k^{m+n}}^{\times}.$ But these mappings are not homomorphisms, so they do not
define the corresponding mappings of classifying spaces.
\end{remark}

So we have to return to fibration (\ref{extended2}).
Let $\mathfrak{bu}_\otimes^i$ be the generalized cohomology theory defined by the
infinite loop space $\BU_{\otimes}$, with the zero term
$\mathfrak{bu}_\otimes^0(X)=[X,\, \BU_{\otimes}]$ \cite{Adams}. Then fibration (\ref{extended2})
corresponds to the exact sequence of Abelian groups
$\mathfrak{bu}_\otimes^0(X)\rightarrow \mathfrak{bu}_\otimes^0(X,\,
\mathbb{Z}[\frac{1}{k}])\rightarrow \mathfrak{bu}_\otimes^0(X;\,
\mathbb{Z}[\frac{1}{k}]/\mathbb{Z})$ associated with the exact sequence of
groups of coefficients $0\rightarrow \mathbb{Z}\rightarrow \mathbb{Z}[\frac{1}{k}]\rightarrow
\mathbb{Z}[\frac{1}{k}]/\mathbb{Z}\rightarrow 0$ (note that $\mathbb{Z}[\frac{1}{k}]/\mathbb{Z}=\varinjlim_n
\mathbb{Z}/k^n\mathbb{Z}$).

Since $\mathfrak{bu}_\otimes^0(X;\, \mathbb{Z}[\frac{1}{k}]/\mathbb{Z})=[X,\, \BN_{k^\infty}^\times]$
is an Abelian group,
we see that the space $\BN_{k^\infty}^\times$ is equipped with the $H$-space structure with
the required property. Put $\FB^k(X):=[X,\, \BN_{k^\infty}^\times]$.
Recall that
the groups $[X,\, \widehat{\Gr}]$ and $[X,\, \BU(k^\infty)]$ we denoted
by $\widehat{\AB}^1(X)$ and $\widehat{\AB}^k(X)$, respectively.
We also defined
the natural transformation (indeed the localization away from $k$) $\widehat{t}_{k*}\colon
\widehat{\AB}^1\rightarrow \widehat{\AB}^k$ induced by $\widehat{t}_k$.
Using the previous arguments we obtain the exact sequence of Abelian groups
$$
\widehat{\AB}^1(X)\stackrel{\widehat{t}_{k*}(X)}{\longrightarrow}\widehat{\AB}^k(X)
\stackrel{\widehat{h}_{k*}(X)}{\longrightarrow}\FB^k(X)
$$
for any finite $CW$-complex $X$.

\begin{definition}
For a finite $CW$-complex $X$ the $k$-{\it primary component} (cf. Remark \ref{kprimcomp})
{\it of
the generalized Brauer group} $\GBr_k(X)$ is the
cokernel $\coker\{ \widehat{h}_{k*}(X)\colon \widehat{\AB}^k(X)
\rightarrow \FB^k(X)\}$.
\end{definition}

Now we want to explain why this definition is a natural generalization of the classical
Brauer group. First of all, let us stress an analogy between the fibration
$$
\mathbb{C}P^\infty \rightarrow \BU(k^\infty)\rightarrow \BPU(k^\infty)
$$
(which relates to the definition of the classical group) and fibration (\ref{extended3})
(or (\ref{extended2})). For example, one can verify that the homotopy sequence of the first
fibration in dimension 2 coincides with the homotopy sequence of the second fibration in all even dimensions.
Therefore in some sense we get a two-periodic generalization of the classical Brauer group.
The role of the Picard group (i.e. the group of isomorphism classes of line bundles
with respect to the tensor product) in our case plays the group of bundles of the form
$(\xi_k,\, (\End(\xi_k),\, \mu,\, \widetilde{M}_{kl})),\; (k,\, l)=1$
(see Subsection 2.6) with respect to the tensor product (cf. Remark \ref{upto1}).

Now it is suitable to consider direct limits (induced by the Kronecker product)
over all natural numbers (not only over powers of
a fixed $k$). In this way we get spaces $\N_\infty^\times:=\varinjlim_n\N_n^\times$ and
$\BN_\infty^\times:=\varinjlim_n\BN_n^\times$ and the fibration (cf. diagram (\ref{usfib}))
$$
\diagram
\BU_{\otimes} \rto & \B{\cal U} \dto &&& \BU_{\otimes}
\rto^{loc\quad \qquad} & \prod_{q\geq 1}\K(\mathbb{Q},\, 2q) \dto \\
& \BN_\infty^\times, & \hbox{i.e.} &&&  \BN_\infty^\times. \\
\enddiagram
$$
Since $\BU_{\otimes}\cong \mathbb{C}P^\infty \times \BSU_{\otimes},$ we see that
this diagram splits:
$$
\diagram
\mathbb{C}P^\infty \dto \rto & \K(\mathbb{Q},\, 2) \rto \dto & \K(\mathbb{Q}/\mathbb{Z},\, 2) \dto \\
\mathbb{C}P^\infty \times \BSU_{\otimes} \rto \dto &
\K(\mathbb{Q},\, 2)\times \prod_{q\geq 2}\K(\mathbb{Q},\, 2q) \rto \dto & \K(\mathbb{Q}/\mathbb{Z},\, 2)
\times \BSN_\infty^\times \dto \\
\BSU_{\otimes} \rto & \prod_{q\geq 2}\K(\mathbb{Q},\, 2q) \rto^{h'} & \BSN_\infty^\times, \\
\enddiagram
$$
where $\BSN_\infty^\times$ is defined by the diagram, $\BN_\infty^\times =\K(\mathbb{Q}/\mathbb{Z},\, 2)
\times \BSN_\infty^\times.$ Notice that the upper fibration correspond to the classical Brauer group,
so the generalized Brauer group is the product of the classical one with the group defined
by the lower fibration.

Since $\B{\cal PU}\simeq \K(\mathbb{Q}/\mathbb{Z},\, 2)\times
\prod_{q\geq 2}\K(\mathbb{Q},\, 2q)$, we obtain the map
$h\colon \B{\cal PU}\rightarrow \BN_{\infty}^\times,\; h=\id_{\K(\mathbb{Q}/\mathbb{Z},\, 2)}\times h'$
(with the homotopy fiber $\Gr =\BSU_{\otimes}$), where $h'$ is defined in the previous diagram,
such that the diagram
\begin{equation}
\label{actgen}
\diagram
\B{\cal U}
\drto^{\widehat{h}} \dto \\
\B{\cal PU}
\rto^{h} &
\BN_{\infty}^\times \\
\enddiagram
\end{equation}
commutes. Put $\overline{\GBr}(X):=\coker\{ h_*(X)\colon \widetilde{\AB}(X)
\rightarrow \FB(X)\}$ (where $\FB(X):=[X,\, \BN_\infty^\times]$).
Then we have the exact sequence of Abelian groups
$$
0\rightarrow \Br(X)\rightarrow
\GBr(X)\rightarrow \overline{\GBr}(X)\rightarrow 0,
$$
where the subgroup $\Br(X)\subset \GBr(X)$
corresponds to those $\N_{k}^\times$-bundles over $X$ whose structure
group can be reduced to $\PU(k)$ (for sufficiently large $k$).

In other words, compared to the
classical case, a whole new step is added to the procedure of the
lifting. At first, we have to
lift an $\N_{k}^\times$-bundle to a $\PU(k)$-bundle.
The second step essentially coincides with the classical case.

\subsection{Application: the twisted $K$-theory}

One of the possible application of the generalized Brauer group is to the twisted $K$-theory.
More precisely, M. Atiyah and G. Segal posed in \cite{Atiyah} the problem of finding of a geometrical
approach to more general twistings than the one comes from projective unitary group $\PU(H)$
of the Hilbert space $H$. Recall that $\PU(H)\simeq \mathbb{C}P^\infty,$ so $\BPU(H)\simeq \K(\mathbb{Z},\, 3).$
For any element $\alpha \in H^3(X; \mathbb{Z})=[X,\, \K(\mathbb{Z},\, 3)]$
(i.e. actually for a $\mathbb{C}P^\infty$-bundle
with the structure group $\PU(H)$),
in particular, for elements of $\Br(X)=H^3_{tors}(X; \mathbb{Z})$,
one can define a ``twisted $K$-group'' $K_{\alpha}^0(X)$ as the group of homotopy classes of sections
of the $\Fred(H)$-bundle associated with $\alpha$ \cite{Atiyah}.
So let us give the analogous construction for bundles with the structure
group $\N_k^\times$.

\begin{remark}
Recall that $\N_k^\times$-bundles represent elements of the generalized Brauer group.
It seems reliable that the twisted $K$-theory associated with such a bundle depends only on the corresponding
element of the generalized Brauer group, i.e. on the class of such a bundle modulo
$\im (\widehat{h}_*)$ (cf. Proposition \ref{lastpropp} below).
\end{remark}

Of course, the group $\N_k^\times \subset M_k({\cal B}(H))$
naturally acts on $H^k$, but it is more fruitful to treat it as an
extension of $\N_1^\times =\N_{k^0}^\times ={\cal B}(H).$ More
precisely, put $\Fred_k:=\N_k\cap
\Fred(H^k)=\pi_k^{-1}(\GL_1(\Delta_k)),$ where $\Fred(H^k)\subset
M_k({\cal B}(H))$ denotes the space of Fredholm operators
contained in $M_k({\cal B}(H)).$ We have the fibration $\Fred_k
\rightarrow \GL_1(\Delta_k)$ with fiber an affine space over
$\ker(\pi_k)=M_k({\cal K}(H)).$ In particular, the projection
$\Fred_k \rightarrow \GL_1(\Delta_k)$ is a homotopy equivalence.
Moreover, there are natural inclusions $\Fred_k \rightarrow
\Fred_{k^2}\rightarrow \ldots$ such that the diagram (consisting
of homotopy equivalences)
$$
\diagram
\Fred_k \rto^{\simeq} \dto_{M_{k}({\cal K}(H))} & \Fred_{k^2} \dto^{M_{k^2}({\cal K}(H))} \\
\GL_1(\Delta_k) \rto^= & \GL_1(\Delta_{k^2}) \\
\enddiagram
$$
commutes.

Note also that there is the commutative diagram
$$
\diagram
\N_k^\times \rto^{\subset} \dto_{\GL_{k}({\cal K}(H))} & \Fred_k^0 \dlto^{M_k({\cal K}(H))} \\
\GL_1(\Delta_k)^0,
\enddiagram
$$
where both down-directed arrows are bundles with the marked fibres.

Clearly, the group $\N_k^\times$ naturally acts on $\Fred_k\subset M_k({\cal B}(H))$ by conjugations.
Moreover, this action is compatible with
the fibration $p_k\colon \Fred_k \rightarrow \GL_1(\Delta_k)$ in the sense that the diagram
\begin{equation}
\label{lastd}
\diagram
\N_k^\times \times \Fred_k \rto^{\quad \varphi_k} \dto_{\id \times p_k} & \Fred_k \dto^{p_k} \\
\N_k^\times \times \GL_1(\Delta_k) \rto^{\quad \; \bar{\varphi}_k} & \GL_1(\Delta_k) \\
\enddiagram
\end{equation}
commutes (where $\varphi_k$ is the mentioned action and $\bar{\varphi}_k$ is the natural action
on the group of invertible elements of Calkin algebra). In particular, the projection $p_k$
is $\N_k^\times$-equivariant, i.e. it commutes with the actions.

Thus for a given $\N_{k}^\times$-cocycle on $X$
we can associate the corresponding $\Fred_k$-bundle $\widetilde{\Fred}_k\rightarrow X$ and
define the corresponding twisted $K$-theory as the set of homotopy classes of
sections of this bundle. Since the space $\Fred_k$ is a topological monoid with respect
to the composition of Fredholm operators, and the group action
preserves the composition, we see that this set is actually a group.

\begin{remark}
Note that the action $\N_{k}^\times \times \Fred(H^k)\rightarrow \Fred(H^k)$ can be extended to the action of
the group $\GL_k({\cal B}(H))\supset \N_{k}^\times$ which is contractible. Therefore for any
$\Fred(H^k)$-bundle with the structure group $\N_{k}^\times$ there exists a fiberwise homeomorphism with the
product bundle $X\times \Fred(H^k)$. But the action
$\varphi_k \colon \N_{k}^\times \times \Fred_k \rightarrow \Fred_k$ cannot be extended to the action
of $\GL_k({\cal B}(H))$ (the latter group does not preserve the subspace $\Fred_k\subset \Fred(H^k)$).
That's why we consider the action $\varphi_k$ on $\Fred_k$.
\end{remark}

\begin{proposition}
\label{lastpropp}
Suppose our $\N_{k}^\times$-cocycle over $X$ comes from some $\GL_k({\cal K}(H))$-cocycle
(recall that $\GL_k({\cal K}(H))$ is a subgroup in $\N_{k}^\times$). Then
the corresponding twisted $K$-theory is isomorphic to the untwisted one $K_{\mathbb{C}}(X)$.
\end{proposition}
{\noindent Proof.\;} Recall that
$\GL_k({\cal K}(H))\subset \N_{k}^\times$ is the kernel of the natural homomorphism
$\N_{k}^\times \rightarrow \GL_1(\Delta_k)^0$. The $\N_k^\times$-equivariant projection $p_k
\colon \Fred_k \rightarrow \GL_1(\Delta_k)$ defines the fiberwise map
$$
\diagram
\widetilde{\Fred}_k \dto \rto^{\widetilde{p}_k\quad} & \widetilde{\GL_1(\Delta_k)} \dlto \\
X, \\
\enddiagram
$$
where $\widetilde{\GL_1(\Delta_k)}\rightarrow X$
is the $\GL_1(\Delta_k)$-bundle associated with the same $\N_{k}^\times$-cocycle as
$\widetilde{\Fred}_k$ (using the action $\bar{\varphi}_k$, see diagram (\ref{lastd})).
Since the subgroup $\GL_k({\cal K}(H))\subset \N_{k}^\times$ is contained in the kernel
of the action $\bar{\varphi}_k\colon \N_k^\times \times \GL_1(\Delta_k)\rightarrow \GL_1(\Delta_k)$, we see that
$\widetilde{\GL_1(\Delta_k)}\rightarrow X$ is a trivial fibration. But the sets
(actually groups) of homotopy classes of sections of two fibrations $\widetilde{\Fred}_k
\rightarrow X$ and $\widetilde{\GL_1(\Delta_k)}\rightarrow X$ are naturally isomorphic,
because $\widetilde{p}_k$ is obviously a fiberwise homotopy equivalence.$\quad \square$

\begin{remark}
Let us remark that our generalized Brauer group generalizes the classical Brauer group
(corresponding to the elements of finite order in $H^3(X; \mathbb{Z})$), but not
the infinite order case corresponding to the Dixmier-Douady class
of a principal $\PU(H)$-bundle \cite{Melrose}. In order to obtain a further generalization we have to
extend sequence (\ref{liftanz3}) up to the sequence $\BGL_{k^\infty}({\cal K}(H))
\stackrel{\widehat{h}^\prime_{k^\infty}}{\longrightarrow}\BN_{k^\infty}^\times \stackrel
{\B\widetilde{\pi}_{k^\infty}}{\longrightarrow}\BGL_1(\Delta_{k^\infty})^0$
(cf. the fibration $\BU(k^\infty)\rightarrow \BPU(k^\infty)\rightarrow \K(\mathbb{Z},\, 3)$).
\end{remark}

\vspace{10mm}

\flushright{e-mail: ershov@higeom.math.msu.su}

\begin{thebibliography}{99}
\bibitem{Adams}
{\sc J.F. Adams:}
Infinite Loop Spaces. Princeton, New Jersey, 1978.
\bibitem{Atiyah}
{\sc M. Atiyah, G. Segal:}
Twisted $K$-theory. arXiv preprint, math.KT/0407054.
\bibitem{Prep}
{\sc A.V. Ershov:} Homotopy theory of bundles with fiber
matrix algebra. {\it Journal of Mathematical Sciences (New York)} Vol.123, No.4 (2004), pp. 4198 - 4220.
\bibitem{e1}
{\sc A.V. Ershov:}
Formal group laws over Hopf algebras and their
application to complex cobordism theory.
Preprint 39 (2002), Max-Planck-Institut f\"ur Mathematik.
\bibitem{e2}
{\sc A.V. Ershov:}
Symmetries in complex cobordism theory
related to stable equivalence classes
of SU-bundles. Preprint 70 (2002),
Max-Planck-Institut f\"ur Mathematik.
\bibitem{Fuchs}
{\sc Ph.A. Griffiths, J.W. Morgan:}
Rational Homotopy Theory and Differential Forms. Birkh\"{a}user, 1981.
\bibitem{Grothendieck}
{\sc A. Grothendieck:}
Le groupe de Brauer I.
{\it Sem. Bourbaki} 1964/65, no. 290, 21p.
\bibitem{Melrose}
{\sc V. Mathai, R.B. Melrose, I.M. Singer:}
The index of projective families of elliptic operators. {\it Geometry $\&$ Topology} Vol.9 (2005), 341-373.
\bibitem{Palais}
{\sc R.S. Palais:}
On the homotopy of certain groups of operators. {\it Topology} 3 (1965), 271-279.
\bibitem{Pirce}
{\sc R.S. Pierce:}
Associative Algebras. Springer Verlag, 1982.
\bibitem{Segal}
{\sc G.B. Segal:}
Categories and cohomology theories.
{\it Topology} 13 (1974), 293--312.
\bibitem{Sullivan}
{\sc D. Sullivan:}
Geometric Topology. MIT, Cambridge, Massachusetts, 1970.
\end{thebibliography}
\end{document}